\let\ORIlabel\label
\let\ORIrefstepcounter\refstepcounter
\AddToHook{package/hyperref/before}{
   \let\label\ORIlabel 
   \let\refstepcounter\ORIrefstepcounter}
\documentclass[hidelinks,onefignum,onetabnum,onealgnum]{siamart251216}

\usepackage{graphicx} %
\usepackage{booktabs}
\usepackage{amsfonts}
\usepackage{amsmath}
\usepackage{amssymb}
 \usepackage{array}
\usepackage{mathtools}
\usepackage{colortbl} %
\usepackage{multicol}
\usepackage{enumitem}
\usepackage{bm}
\usepackage{comment}
\usepackage[font=small]{caption}
\usepackage{pgfplots}
\usepackage{placeins}
\definecolor{CeruleanRef}{RGB}{12,127,172}
\hypersetup{colorlinks=true,allcolors=CeruleanRef}
\hypersetup{pdfstartview=FitB,pdfpagemode=UseNone}

\usepackage{rotating}
\RequirePackage{algorithm}[0.1]

\usepackage{algpseudocode}
\algnewcommand{\Initialize}[1]{
  \State \textbf{Initialize:}
  \Statex \hspace*{\algorithmicindent}\parbox[t]{.8\linewidth}{\raggedright #1}
}
\algnewcommand{\Indent}[2]{
  \State {#1}
  \vspace{-2mm}
  \Statex \hspace*{\algorithmicindent}\parbox[t]{.9\linewidth}{\raggedright #2}
}

\usepackage{subfig}
\usepackage{hhline}
\usepackage{multirow}
 \usepackage{wrapfig}
 \usepackage{standalone}
\usepackage{tikz}
\usetikzlibrary{trees}
\usetikzlibrary{calc}

\newcommand{\vectt}[1]{\boldsymbol{\mathbf{#1}}}

\newcommand{\dx}{\mathrm{d}x}
\newcommand{\ds}{\mathrm{d}s}

\newcommand{\supp}{\mathrm{supp}}

\newcommand{\divv}{\operatorname{div}}

\DeclareMathOperator*{\argmin}{arg\,min}

\crefname{problem}{Problem}{Problems}

\title{Preconditioning proximal Galerkin: \\ mesh and degree robust solvers for variational inequalities}

\author{
Ioannis P.~A.~Papadopoulos\thanks{%
Mathematical Institute, University of Oxford, UK}
(\email{ioannis.papadopoulos@maths.ox.ac.uk})
\and
Pablo D.~Brubeck\thanks{%
Mathematical Institute, University of Oxford, UK}
(\email{pablo.brubeck@maths.ox.ac.uk})
}

\begin{document}

\maketitle

\begin{abstract}
The proximal Galerkin algorithm is a structure-preserving solver for variational inequalities that results in discretization-independent nonlinear iteration counts. However, in large three-dimensional simulations, the arising Newton saddle point systems cannot be feasibly solved by direct methods. To leverage an iterative linear solver, we propose a robust and scalable block preconditioner.
Instead of using the usual dual Schur complement that acts on the latent variable, we design a preconditioner for the primal Schur complement. This leads to a block preconditioner that requires two approximate inverses of the latent block and only one of the primal Schur complement. The latent block represents an algebraic equation which, after a small perturbation, can be cheaply inverted. We then introduce two choices of approximations for the primal Schur complement: a PDE-operator-based and an algebraic approach. Either approximation can be efficiently inverted by geometric multigrid. 

We investigate the effectiveness of these strategies when proximal Galerkin is applied to obstacle, Signorini, and strain-constrained elasticity problems and also compare with a monolithic multigrid approach. Choices of standard relaxation methods, in the block preconditioner route, lead to robust iteration counts with respect to the mesh size and discretization degree. Moreover, in the context of the Signorini, specialized relaxation and prolongation also result in effective solvers for nearly-incompressible elastic materials. The linear solves are parameter-robust in 2D, with only mild growth in 3D. This allows for the solution of a three-dimensional contact problem, with a Poisson's ratio of 0.499, discretized with over ten million degrees of freedom.

\end{abstract}



\section{Introduction}
\label{sec:introduction}

The latent variable proximal point (LVPP) method is a framework for solving variational inequalities (VIs) such as, but not limited to, obstacle, Signorini contact, and elastic-plastic torsion problems \cite{dokken2024,kinderlehrer2000}. LVPP is a saddle point reformulation of the infinite-dimensional variant of the Bregman proximal point (BPP) algorithm where the reformulation is achieved by introducing a latent variable denoted by $\psi$. LVPP and BPP are equivalent on the infinite-dimensional level and result in smooth (systems of) PDEs as the first-order optimality conditions to each subproblem. However, after discretization, the saddle point nonlinear PDEs of LVPP require significantly fewer iterations from Newton's method to reach convergence. As such, LVPP can be thought of as a choice of lifted Newton's method.

LVPP offers significant flexibility with the choice of discretization.  When coupled with a Galerkin discretization, such as a finite element method (FEM), the resulting algorithm is referred to as proximal Galerkin (PG) \cite{dokken2024,ern2026,fu2026,fu2024,keith2023,papadopoulos2024}. Theoretical and experimental evidence shows that PG delivers robust convergence with mesh- and degree-independent outer Newton iteration counts \cite{dokken2024,keith2025}. This paper introduces a multigrid-based block preconditioner for the iterative linear solvers applied to each Newton system \cite{benzi2005,wathen2015}. Experimentally, we observe the proposed iterative solvers deliver mesh- and degree-robust Krylov iteration counts. With specialized smoothing and prolongation, we can solve nearly-incompressible Signorini problems with parameter-robust Krylov iteration counts in 2D, with only mild growth in 3D \cite{farrell2021b}.

In summary, we propose an algorithm for solving VIs  where the number of Newton iterations and inner Krylov iterations (to solve the Newton linear systems) are mesh- and degree-robust. 

In \Cref{sec:pg} we introduce the PG algorithm for the obstacle \cite{stampacchia1964}, Signorini \cite{fichera1964}, and strain-constrained elasticity problems \cite[Ex.~7]{rodrigues2019}. The Newton saddle point systems are derived in \Cref{sec:discretization} where their FEM discretization is introduced. \Cref{sec:preconditioning} then details the preconditioners where the block preconditioner is introduced \Cref{sec:block}. Subspace correction methods are then described in \Cref{sec:domain-decomp} and a monolithic multigrid strategy in \Cref{sec:monolithic}. We  investigate the effectiveness of the preconditioned PG solver via numerical examples in \Cref{sec:examples} and give our conclusions in \Cref{sec:conclusions}.

\section{Proximal Galerkin}
\label{sec:pg}

In this work we consider open, bounded, and Lipschitz domains $\Omega \subset \mathbb{R}^d$ for $d \in \{2,3\}$. We denote the usual Lebesgue spaces by $L^p(\Omega)$, $p \in [1,\infty]$, and the Sobolev spaces by $W^{s,p}(\Omega)$, $s \geq0$, $p \in [1,\infty]$ \cite{Adams2003}. We define the Hilbert spaces $H^s(\Omega) \coloneqq W^{s, 2}(\Omega)$. For a non-empty $\Gamma_D \subseteq \partial \Omega$, we denote  $H^1_{\Gamma_D,g}(\Omega) \coloneqq \{ u \in H^1(\Omega) : u|_{\Gamma_D} = g\}$ where $g \in H^1(\Gamma_D)$ and $|_{\Gamma_D} : W^{1,p}(\Omega) \to W^{1-1/p,p}(\Gamma_D)$ denotes the trace operator.  Moreover, $H^1_0(\Omega) = H^1_{\partial \Omega, 0}$. We use $H^s(\Omega)^d$ for a vector-valued Sobolev space with $d$ copies and $L^p(\Omega)^{d \times d}_{\mathrm{sym}}$, $p \geq 1$, to denote the $d \times d$ symmetric tensor-valued Lebesgue space. The inner product $(\cdot,\cdot)$ without a subscript always denotes the $L^2(\Omega)$-inner product. 

Let $V$ denote a Hilbert space and consider the following pointwise constrained optimization problem: find $u \in V \cap K$ such that
\begin{align}
u =\argmin_{v \in V \cap K} j(v).
\label{eq:opt}
\end{align}
Here $j$ is a Fr\'echet differentiable functional $j : V \to \mathbb{R}$ and $K$ is a closed convex set. We consider three choices for $K$ in this work:
\begin{subequations}
\begin{align}
K_o &\coloneqq \{u \leq \phi \; \text{a.e.~in} \; \Omega\},\label{K:obstacle}\\
K_s &\coloneqq  \{u|_{\Gamma_C} \cdot n \leq \phi \ \; \text{a.e.~on} \; \Gamma_C\}, \label{K:signorini}\\
K_g &\coloneqq \{|\varepsilon(u)| \leq \phi \; \text{a.e.~in} \; \Omega\}. \label{K:strain}
\end{align}
\end{subequations}
We refer to these constraints as obstacle, Signorini, and strain constraints, respectively. $\phi$ denotes a scalar-valued function, $\varepsilon = (\nabla + \nabla^\top)/2$ is the symmetric gradient operator, $\Gamma_C \subset \partial \Omega$ is an open, nonempty, connected part of the domain boundary, $n \in L^\infty(\Gamma_C)^d$ is a problem-given unit-vector function and $u : \Omega \to \mathbb{R}$ is scalar-valued when considering obstacle constraints and vector-valued $u : \Omega \to \mathbb{R}^d$ otherwise.

\subsection{PG for the obstacle problem} PG was first introduced for obstacle constraints in \cite{keith2023}. Let $V = H^1_0(\Omega)$ and $\phi \in H^1(\Omega)$ with $\phi|_{\partial \Omega} \geq 0$.  The \emph{obstacle problem} seeks $u \in V \cap K_o$ such that
\begin{align}
u = \argmin_{v \in V \cap K_o}  j_1(v) \coloneqq \int_\Omega \frac{1}{2}|\nabla v|^2 - f v \, \dx,
\label{eq:obstacle}
\end{align} 
for a given $f \in L^2(\Omega)$.

Given a so-called latent variable $\psi^{k-1} \in L^\infty(\Omega)$, for $k \in \mathbb{N}$, the PG algorithm iteratively solves the following nonlinear saddle point problem to find $(u^k, \psi^k) \in H^1_0(\Omega) \times L^\infty(\Omega)$ that, for all $(v, q) \in H^1_0(\Omega) \times L^\infty(\Omega)$, satisfies
\begin{align}
\begin{split}
\alpha_k (\nabla u^k, \nabla v) + (\psi^k, v) &= \alpha_k (f, v) + (\psi^{k-1}, v)\\
(u^k + \exp(-\psi^k), q) &= (\phi, q).
\end{split} \label{eq:pg:obstacle}
\end{align}
We refer to each nonlinear solve of \cref{eq:pg:obstacle} as a \emph{proximal step}. If $u^* \in H^1_0(\Omega) $ denotes the unique solution to the obstacle problem, then the iterates $u^k$ satisfy the error bound $\|\nabla(u^k - u^*)\|^2_{L^2(\Omega)} \leq c \left(\sum_{j=1}^k \alpha_j \right)^{-1} \|\nabla(u^0 - u^*)\|^2_{L^2(\Omega)}$. Hence, the algorithm converges provided $\sum_{j=1}^k \alpha_j \to \infty$. Note that the sum diverges, and thus the algorithm converges, even if $\alpha_k = 1$ for all $k \in \mathbb{N}$. The faster $\alpha_k$ grows, the quicker the theoretical convergence to the solution. In practice, \cref{eq:pg:obstacle} is solved via Newton's method. Large growth rates of $\alpha_k$ may lead to faster outermost convergence but at the cost of making each proximal step (solve of \cref{eq:pg:obstacle}) harder to solve. This manifests as requiring more Newton steps to converge. Hence, one strikes a balance between the rate of proximal convergence and the difficulty of solving the proximal subproblems. 

\subsection{PG for the Signorini problem}

In the \emph{Signorini problem}, $u \in V = H^1_{\Gamma_D, g}(\Omega)^d$, $d \in \{2,3\},$ is a vector-valued function representing the deformation of an elastic body from its original position and $g \in H^1(\Gamma_D)^d$ is a prescribed Dirichlet boundary condition. Let $\Gamma_C \subseteq \partial \Omega$ denote the candidate region for contact on the boundary of the body, $n \in L^\infty(\Gamma_C)^d$ be a given outward unit-norm direction defined on $\Gamma_C$, and $\phi \in L^\infty(\Gamma_C)$ be a prescribed gap function. Then,  given a body force $f \in L^2(\Omega)^d$, the Signorini problem seeks $u \in V \cap K_s$ satisfying
\begin{align}
u = \argmin_{v \in V \cap K_s} j_2(v) \coloneqq \int_\Omega \mu |\varepsilon(v)|^2 + \frac{\lambda}{2} |\divv v |^2 - f \cdot v \, \dx.
\label{eq:signorini}
\end{align}
In this context, the symmetric gradient $\varepsilon(u)$ models the infinitesimal strain of the material and $\mu, \lambda$ are referred to as the Lam\'e parameters. For a fixed $\mu>0$, the limit $\lambda \to \infty$ corresponds to an increasingly incompressible material.

PG, for the Signorini problem, iteratively solves the following nonlinear saddle point problem: find $(u^k, \psi^k) \in H^1_{\Gamma_D,g}(\Omega)^d \times L^\infty(\Gamma_C)$ such that for all $(v, q) \in H^1_{\Gamma_D,0}(\Omega)^d \times L^\infty(\Gamma_C)$ \cite[Sec.~3.2]{dokken2024}
\begin{align}
\begin{split}
2\mu\alpha_k (\varepsilon(u^k), \varepsilon(v)) +  \lambda \alpha_k(\divv u^k, \divv v) + \int_{\Gamma_C} (\psi^k-\psi^{k-1}) (v \cdot n) \, \ds &= \alpha_k (f, v),  \\
\int_{\Gamma_C} (u^k \cdot n+ \exp(-\psi^k) )q \, \ds & = \int_{\Gamma_C} \phi q \, \ds.
\end{split}  \label{eq:pg:signorini}
\end{align}
The latent variable $\psi^k$ is a scalar-valued function with support on a subdomain of codimension one. Thus its discretization typically requires significantly fewer degrees of freedom (DoFs) than $u$.

\subsection{PG for strain-constrained elasticity}
Let $V = H^1_{\Gamma_D,g}(\Omega)^d$, $d \in \{2,3\}$, where $g \in H^1(\Gamma_D)^d$ is a prescribed Dirichlet boundary condition. Given a strain-cap $\phi \in L^\infty(\Omega)$, we seek to find $u \in V \cap K_g$ that satisfies
\begin{align}
u = \argmin_{v \in V \cap K_g} j_2(v).
\label{eq:strain}
\end{align}
This problem models elastic materials that undergo a locking phenomenon where, after a certain threshold,  the strain of a material remains fixed despite increases in the stress \cite[Ex.~7]{rodrigues2019}.

In each proximal step we seek $(u^k,\psi^k) \in H^1_{\Gamma_D,g}(\Omega)^d \times L^1(\Omega)^{d \times d}_{\mathrm{sym}}$  such that for all $(v, q) \in H^1_{\Gamma_D,0}(\Omega)^d \times L^1(\Omega)^{d \times d}_{\mathrm{sym}}$
\begin{align}
\begin{split}
2\mu\alpha_k (\varepsilon(u^k), \varepsilon(v)) +  \lambda \alpha_k(\divv u^k, \divv v) + (\psi^k, \varepsilon(v)) &= \alpha_k (f, v) + (\psi^{k-1}, \varepsilon(v))\\
\gamma(\varepsilon(u^k) - R_\phi(\psi^k), q) &= 0,
\end{split} \label{eq:pg:strain}
\end{align}
where
\begin{align}
R_\phi(\psi) = \frac{\phi \psi}{\sqrt{1+|\psi|^2}}.
\end{align}
Here $\gamma>0$ is a constant that helps better enforce the constraints after discretization. In exact arithmetic, it leaves the solution unchanged. In the experiments we consider values of $\gamma=100$ in 2D and $\gamma=5000$ in 3D. 

In direct contrast to the previous two problems, $\psi$ is now a tensor-valued function. Theoretical guarantees of convergence of PG for problems with gradient-type constraints is an open problem although experimental evidence can be found in \cite{dokken2024,papadopoulos2024}.

\section{Newton linearization and FEM discretization}
\label{sec:discretization}
Each proximal step requires the solve of a nonlinear system of PDEs. In this work, we solve each nonlinear system with Newton's method, resulting in a sequence of linear systems. More precisely, we first perform a Newton linearization of \cref{eq:pg:obstacle,eq:pg:signorini,eq:pg:strain} and discretize the resulting infinite-dimensional linear system with an FEM. 

\begin{figure}[h!]
\centering
\includegraphics[width =0.6 \textwidth]{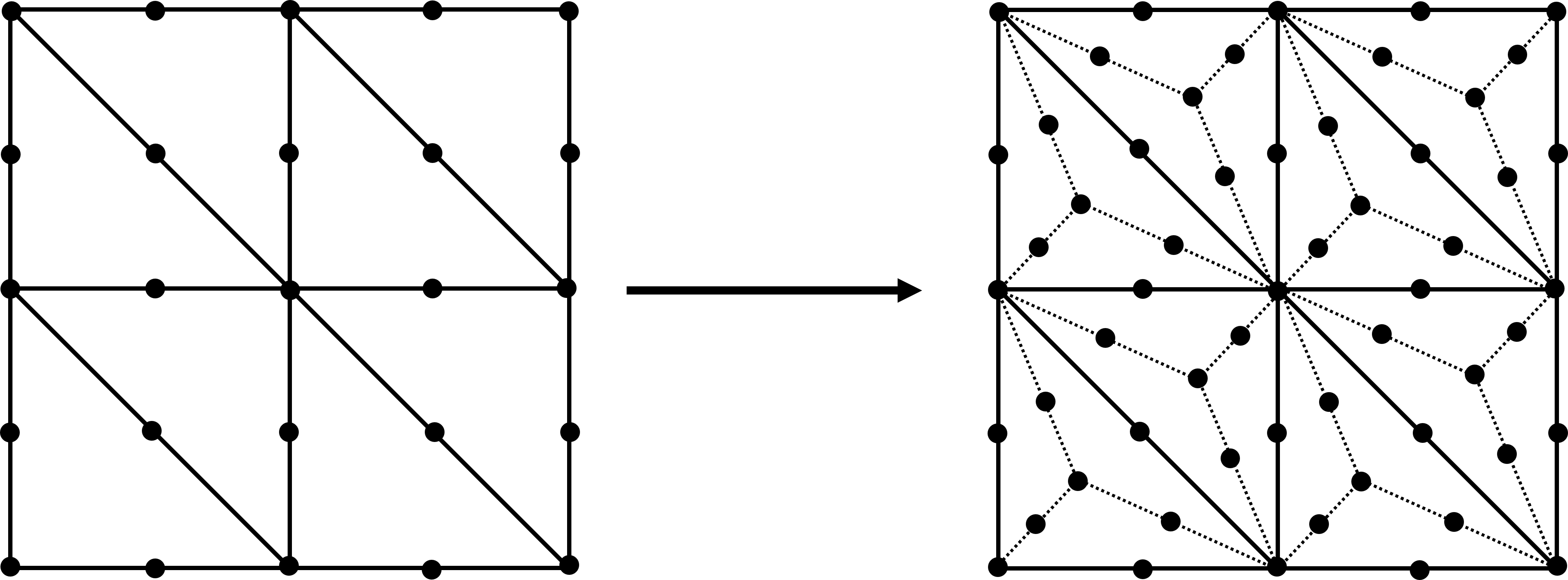} 
\caption{The DoFs of a $\mathrm{CG}_2$ (left) and $\mathrm{CG}_2$-A (right) discretization.}
\label{fig:alfeld}
\end{figure}

The domain $\Omega$ is triangulated by a quasi-uniform mesh of simplicial cells, denoted $\mathcal{T}_h$. We shall consider a variety of FEM spaces for $u^k \in V_h$ and $\psi^k \in Q_h$ including a continuous Lagrange element of degree $p$ ($\mathrm{CG}_p$), discontinuous Lagrange ($\mathrm{DG}_p$) \cite{defelement}, and Mardal--Tai--Winther elements $(\mathrm{MTW}_1)$ \cite{mardal2002,defelement}. We also consider discretizations where each cell in the mesh has undergone an Alfeld split \cite{brubeck2026,farrell2022} as demonstrated in \cref{fig:alfeld}. These are denoted $\mathrm{CG}_p$-A, $p \geq d$, for a degree $p$ continuous Lagrange element over an Alfeld split.

Let $u^k_h \in V_h$ and $\psi^k_h \in Q_h$ denote the finite element approximations of $u^k$ and $\psi^k$ and $\vectt{u} \in \mathbb{R}^{N_u}$ and $\vectt{\psi}\in \mathbb{R}^{N_\psi}$ denote the corresponding FEM coefficient vectors (where we drop the superscript $^k$ for brevity). Then the resulting discretized Newton linear systems are of the saddle point form:
\begin{align}
J\vectt{\delta} = \begin{pmatrix}
A_\alpha  & B\\
B^\top & -D
\end{pmatrix}
\begin{pmatrix}
\vectt{\delta}_u\\
\vectt{\delta}_\psi
\end{pmatrix}= \vectt{b}.
\label{eq:saddle-point}
\end{align}
$\vectt{\delta}_u \in \mathbb{R}^{N_u}$ and $\vectt{\delta}_\psi \in \mathbb{R}^{N_\psi}$ denote the Newton updates for $\vectt{u}$ and $\vectt{\psi}$, respectively, and $\vectt{b}$ denotes the discretization of the negative residual. $A_\alpha$ is a symmetric positive-definite (SPD) matrix corresponding to the $\alpha$-scaled stiffness matrix associated with the problem, $B$ and $B^\top$ characterize the interaction between the $u$ and $\psi$ DoFs (identity, trace, and $\varepsilon$/$\divv$ interactions for \cref{eq:pg:obstacle,eq:pg:signorini,eq:pg:strain}, respectively) and $D$ is an SPD matrix which captures the linearization of the exponential or $R_\phi$-terms. We refer to $J$ in \cref{eq:saddle-point} as the \emph{Jacobian}. For \cref{eq:pg:obstacle,eq:pg:signorini,eq:pg:strain}, $A_\alpha$ discretizes a differential operator, while $D$ discretizes an algebraic operator.

The precise entries are:
\begin{subequations}
\begin{align}
[A_\alpha]_{ij} 
&= \begin{cases}
 \alpha_k \int_\Omega \nabla \varphi_i \cdot \nabla \varphi_j \, \dx & \text{Obstacle \cref{eq:pg:obstacle}},\\
 \alpha_k \int_\Omega 2\mu \varepsilon(\varphi_i) : \varepsilon(\varphi_j) +  \lambda \divv \varphi_i \cdot \divv \varphi_j \, \dx  & \text{Signorini \cref{eq:pg:signorini}},\\
 \alpha_k \int_\Omega 2\mu \varepsilon(\varphi_i) : \varepsilon(\varphi_j) +  \lambda \divv \varphi_i \cdot \divv \varphi_j \, \dx  & \text{Strain \cref{eq:pg:strain}},
\end{cases}\\
B_{ij} &= \begin{cases}
\int_\Omega \varphi_i \eta_j \, \dx  & \hspace*{40mm} \text{Obstacle \cref{eq:pg:obstacle}},\\
  \int_{\Gamma_C} (\varphi_i \cdot n) \eta_j \, \ds & \hspace*{40mm}\text{Signorini \cref{eq:pg:signorini}},\\
 \int_\Omega \varepsilon(\varphi_i) : \eta_j \, \dx  & \hspace*{40mm}\text{Strain \cref{eq:pg:strain}},
\end{cases}\\
D_{ij} &= \begin{cases}
\int_\Omega \exp(-\psi^{k, \ell}_h) \eta_i \eta_j \, \dx & \text{Obstacle \cref{eq:pg:obstacle}},\\
\int_{\Gamma_C} \exp(-\psi^{k, \ell}_h) \eta_i \eta_j \, \ds   & \text{Signorini \cref{eq:pg:signorini}},\\
 \int_\Omega \frac{\phi}{s^{k,\ell}_h} \eta_i : \eta_j - \frac{\phi}{(s^{k,\ell}_h)^{3}} (\psi^{k, \ell}_h : \eta_i) (\psi^{k, \ell}_h : \eta_j) \, \dx & \text{Strain \cref{eq:pg:strain}},
\end{cases} \label{eq:D}
\end{align}
\end{subequations}
where $\varphi_i$, $i=1,\dots,N_u$ and $\eta_j$, $j=1,\dots,N_\psi$ denote the FEM basis functions for $u$ and $\psi$, respectively, $\psi^{k,\ell}_h$ is the $\ell$-th Newton iterate for $\psi^k_h$ in the $k^\text{th}$ proximal step, and $s^{k,\ell}_h = (1+|\psi^{k,\ell}_h|^2)^{1/2}$.

\begin{figure}[ht]
\centering
 \begin{tikzpicture}[%
every node/.style={draw=black, thick, anchor=west},
grow via three points={one child at (-1.5,-0.7) and
    two children at (-1.5,-0.7) and (-1.5,-1.4)},
edge from parent path={(\tikzparentnode.190) |- (\tikzchildnode.west)}]
\node {Constrained Optimization Problem \cref{eq:opt}}
child { node {System of Nonlinear PDEs (Proximal Galerkin) \cref{eq:pg:obstacle,eq:pg:signorini,eq:pg:strain}}
    child { node {Saddle Point Linear System  \cref{eq:saddle-point} (Newton \& Discretization)}}
};
\end{tikzpicture}
\caption{PG strategy for solving the constrained optimization problem \cref{eq:opt}.}\label{fig:pg}
\end{figure}
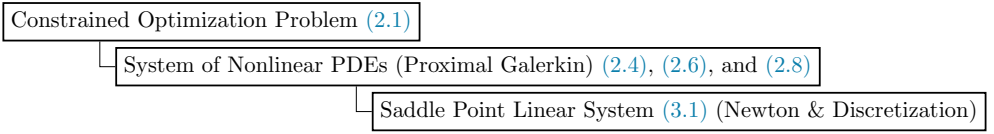

In \cref{fig:pg}, we summarize the PG strategy for solving \cref{eq:opt}.

We subsequently wish to solve the Newton system \cref{eq:saddle-point} via an iterative Krylov method, namely flexible GMRES (FGMRES) \cite{saad1993}. In order for FGMRES to converge within a reasonable number of iterations, we require a suitable easy-to-invert preconditioner $J_P$ for $J$ in \cref{eq:saddle-point}.

\section{Preconditioning}
\label{sec:preconditioning}
In this section we introduce the block preconditioning strategy for the FGMRES method when solving the PG Newton systems \cref{eq:saddle-point}  in \Cref{sec:block}. We will leverage standard iterative  methods for the multigrid relaxation such as Jacobi \cite[Ch.~4.1]{saad2003}, block-Jacobi \cite[Alg.~4.1]{saad2003}, Chebyshev \cite[Alg.~12.1]{saad2003}, preconditioned conjugate gradient (PCG) \cite[Ch.~6.7]{saad2003}, and GMRES \cite[Ch.~6.5]{saad2003}. However, we will also investigate iterative solvers from the subspace correction literature such as the vertex-star patch and Vanka-type iteration. Since these are less commonly known, we introduce them in \Cref{sec:domain-decomp}. We also propose a monolithic strategy for comparison purposes in \Cref{sec:monolithic}.

\subsection{Block preconditioning}
\label{sec:block}
Block preconditioners are a common strategy for solving saddle point systems \cite{benzi2005,pearson2017,rees2010,rees2021,wathen2015}. There are two choices of block preconditioner based on either an LDU ($J_D$) or UDL decomposition $(J_P)$  \cite[Sec.~5]{benzi2005}:
\begin{align}
J_D^{-1} &=
\begin{pmatrix}
I & -\hat{A}_\alpha^{-1} B \\ 0& I
\end{pmatrix}
\begin{pmatrix}
\hat{A}_\alpha^{-1}& 0 \\ 0 &  \hat{S}_D^{-1}
\end{pmatrix}
\begin{pmatrix}
I &  0 \\ -B^\top \hat{A}_\alpha^{-1} & I
\end{pmatrix},\\
J_P^{-1} &=
\begin{pmatrix}
I & 0 \\  \hat{D}^{-1} B^\top & I
\end{pmatrix}
\begin{pmatrix}
\hat{S}_P^{-1}& 0 \\ 0 &  -\hat{D}^{-1}
\end{pmatrix}
\begin{pmatrix}
I &  B \hat{D}^{-1} \\  0 & I
\end{pmatrix},
\label{eq:schur2}
\end{align}
where $\hat{A}_\alpha^{-1}$, $\hat{D}^{-1}$, $\hat{S}_D^{-1}$, and $\hat{S}_P^{-1}$ denote approximations to the actions of $A_\alpha^{-1}$, $D^{-1}$, $S_D^{-1}\coloneqq (-D - B^\top A_\alpha^{-1} B)^{-1}$ and $S^{-1}_P \coloneqq (A_\alpha + B D^{-1} B^\top)^{-1}$, respectively. $S_D$ and $S_P$ are called the \emph{dual} Schur complement and \emph{primal} Schur complement, respectively.

Each application of $J_P^{-1}$ requires two actions of $\hat{D}^{-1}$ and one of $\hat{S}_P^{-1}$ whereas each application of $J_D^{-1}$ requires two actions of $\hat{A}_\alpha^{-1}$ and one of $\hat{S}_D^{-1}$. If the approximations are exact, i.e.~$\hat{A}^{-1}_\alpha = A^{-1}_\alpha$, $\hat{D}^{-1} = D^{-1}$, $\hat{S}^{-1}_D = S^{-1}_D$, and $\hat{S}^{-1}_P = S^{-1}_P$, then
\[
J_D^{-1} J = J_P^{-1} J = I.
\]

The first block preconditioner for PG leveraged the LDU decomposition $J_D$ \cite{papadopoulos2024}. In the $hp$-framework of \cite{papadopoulos2024}, the size of $A_\alpha$ typically remained small and could be handled by a Cholesky factorization. The dual Schur complement $\hat{S}_D = S_D$ was inverted via GMRES (its action could be applied efficiently) and a preconditioner for $S_D$ was constructed by dropping the lowest-order basis functions in the discretization of $u$ resulting in a block-diagonal preconditioner. For a sufficiently small GMRES tolerance, a single action of $J_D^{-1}$ was sufficiently accurate to approximate $J^{-1}$.

This strategy resulted in robust iteration counts with respect to $p$-refinement but was less robust to $h$-refinement. Each GMRES iteration for $S_D$ requires the action of $A_\alpha^{-1}$.  When a factorization of $A_\alpha$ is unavailable and $A_\alpha^{-1}$ is computed iteratively, the nested GMRES iteration makes this solver computationally intractable. A further issue is that it is unclear how to construct a $h$-robust preconditioner for $S_D$ since, for the obstacle problem, $S_D$ is a finite-dimensional approximation of the nonlocal operator $S_D \approx -\exp(-\psi^{k,\ell}) -\alpha^{-1}_k (-\Delta)^{-1}$ where the superscripts denote the $\ell^\text{th}$ Newton iterate of the $k^\text{th}$ proximal step. These observations motivate a different strategy here.

In this work, we instead consider the UDL decomposition $J_P$. Since $D$ is a sparse and SPD discretization of an algebraic operator, a seemingly natural choice for $\hat{D}$ is $\hat{D} = D$. In the first few proximal steps this is effective. However, in the later proximal steps, one encounters regions in $\Omega$ where $\psi^{k,\ell}_h \gg 1$. Examining \cref{eq:D}, this leads to entries in $D$ becoming near zero and, in turn, $D$ becomes increasingly singular. This is disastrous for the effectiveness of $J_P$. Hence, we instead stabilize the action of $\hat{D}^{-1}$ by introducing a small perturbation 
\begin{align}
\hat{D} = D + \epsilon M, \quad \text{for $0 < \epsilon \ll 1$,} \label{eq:D-hat}
\end{align}
where $M$ is the $\psi$-mass matrix $M_{ij} = (\eta_i, \eta_j)$ for the obstacle and strain-constrained problems and $\int_{\Gamma_C} \eta_i \eta_j \ds$ for the Signorini problem. 
A key point is that the true Jacobian $J$ is left unperturbed (i.e.~only the latent block in the \emph{preconditioner} $J_P$ is perturbed), and hence we are still leveraging a true Newton method (rather than a quasi-Newton method). Examining $\hat{D}$, we notice it is effectively a weighted mass matrix. As a consequence it can be approximately inverted efficiently by a handful of PCG iterations preconditioned with a Jacobi or block-Jacobi iteration \cite{wathen1987}.

We now consider two approximations $\hat{S}_P$ for the Schur complement $S_P$. 

\subsubsection*{Algebraic Schur complement approximant} If $\psi$ is discretized with a $\mathrm{DG}$ basis, then $\hat{D}$ is block-diagonal, its inverse can be computed cell-wise, and the triple-matrix product $B \hat{D}^{-1} B^\top$ can be assembled exactly and cheaply. If instead a continuous basis is used for $\psi$, $\hat{D}^{-1}$ is dense and $S_P$ cannot be formed directly. We therefore introduce the approximations $\hat{D}_a$ and $B_a$ of $\hat{D}$ and $B$, obtained by assembling these operators with respect to the broken (discontinuous) counterpart of the latent finite element space, and define the algebraic Schur complement approximant
\begin{align}
\label{eq:Sd}
\hat{S}_{P,a} \coloneqq A_\alpha + B_a \hat{D}_a^{-1} B_a^\top = A_\alpha + \sum_{K \in \mathcal{T}_h} L_K^\top B_K \hat{D}_K^{-1} B_K^\top L_K,
\end{align}
where $B_K$ and $\hat{D}_K$ denote the element matrices on the cell $K$ and $L_K$ is the local-to-global map for the primal degrees of freedom. In other words, $\hat{S}_{P,a}$ is an element-based Schur complement in the spirit of Rees and Wathen \cite{rees2021}: a local Schur complement is formed on each cell and the contributions are summed into the global primal space. In practice, $B_a$ and $\hat{D}_a$ are never assembled as global matrices; the dense element inverses and products are evaluated cell-by-cell inside a single assembly kernel generated by Firedrake's Slate language \cite{slate}. Consequently, $\hat{S}_{P,a}$ has the same sparsity pattern as $A_\alpha$, costs roughly as much to construct as an ordinary finite element matrix, and is amenable to the same geometric multigrid techniques as $A_\alpha$. 

\subsubsection*{Operator Schur complement approximant}  An alternative approach is inspired by operator preconditioning \cite{hiptmair2006,mardal2011}. In all three problems, $\hat{D}$ represents an algebraic operator containing no differential operators. Hence, on the continuous level, the action of $\hat{D}^{-1}$ is equivalent to a pointwise inverse. For example, in the obstacle and Signorini problems:
\begin{align*}
\hat{D} \approx \epsilon + \exp(-\psi^{k,\ell})  \implies \hat{D}^{-1} \approx \frac{1}{\epsilon + \exp(-\psi^{k,\ell})}.
\end{align*}
For the strain-constrained problem:
\begin{align*}
&\hat{D} \approx \left(\epsilon + \frac{\phi }{s}\right) I \otimes I -\frac{\phi}{s^3}  (\psi^{k, \ell}_h  \otimes \psi^{k, \ell}_h) \\ &\implies 
\hat{D}^{-1} \approx \frac{s^{k,\ell}_h}{\epsilon s^{k,\ell}_h + \phi} \left( I \otimes I + \frac{\phi}{\epsilon (s^{k,\ell}_h)^3 + \phi} (\psi^{k, \ell}_h  \otimes \psi^{k, \ell}_h)  \right),
\end{align*}
where $I \in \mathbb{R}^{d \times d}$ denotes the identity operator and $s^{k,\ell}_h = (1+|\psi^{k,\ell}_h|^2)^{1/2}$.  Here, the pointwise inverse of the tensor-valued operator is derived by noticing that $\hat{D}$ is a rank-one perturbation of a scaled identity operator and applying the Sherman--Morrison formula.

In turn, this allows us to approximate $S_P$ by discretizations of a PDE operator.  For the obstacle, and Signorini problems, respectively, the operator Schur complement approximant, denoted $\hat{S}_{P,o}$, takes the entries:
\begin{subequations}
\begin{align}
[\hat{S}_{P,o}]_{ij} &= \int_\Omega \alpha_k \nabla \varphi_i \cdot \nabla \varphi_j + \frac{\varphi_i \varphi_j}{\epsilon + \exp(-\psi^{k,\ell}_h)} \, \dx,\\
\begin{split}
[\hat{S}_{P,o}]_{ij} &= \alpha_k\int_\Omega  2 \mu \varepsilon( \varphi_i) : \varepsilon( \varphi_j) + \lambda \divv \varphi_i \cdot \divv \varphi_j \, \dx \\
&+ \int_{\Gamma_C}\frac{(\varphi_i \cdot n) (\varphi_j \cdot n)}{\epsilon + \exp(-\psi^{k,\ell}_h)} \, \ds,
\end{split}
\end{align}
whereas for strain-constrained problems,
\begin{align}
\begin{split}
[\hat{S}_{P,o}]_{ij} &=  \alpha_k\int_\Omega  2 \mu \varepsilon( \varphi_i) : \varepsilon( \varphi_j) + \lambda \divv \varphi_i \cdot \divv \varphi_j \, \dx\\
& +  \int_\Omega\frac{s^{k,\ell}_h}{\epsilon s^{k,\ell}_h+\phi}\left( \varepsilon(\varphi_i) : \varepsilon(\varphi_j) \right)\dx \\
& + \int_\Omega \frac{\phi s^{k,\ell}_h}{(\epsilon s^{k,\ell}_h+\phi)(\epsilon (s^{k,\ell}_h)^3+\phi)} (\varepsilon(\varphi_i) : \psi^{k,\ell}_h) (\varepsilon(\varphi_j) : \psi^{k,\ell}_h) \, \dx,
\end{split}
\end{align} \label{eq:So}
\end{subequations}
where $\varphi_i$, $i=1,\dots,N_u$ are the FEM basis functions for $u$.

Both $\hat{S}_{P,a}$ and $\hat{S}_{P,o}$ are SPD matrices representing a second-order linear PDE with variable coefficients. Hence, we approximate the actions of $\hat{S}^{-1}_{P,a}$ and $\hat{S}^{-1}_{P,o}$ by geometric multigrid. On each coarser level, the problems are rediscretized to construct the coarser analogues. Effective multigrid cycles require adequate relaxation on the finer grids. In \Cref{sec:examples}, we shall observe that $h$ and $p$-robust iteration counts are attained by applying standard relaxation methods (Chebyshev or PCG).  Moreover, parameter-robust iteration counts can be achieved via vertex-star patch iterations which we now introduce in \Cref{sec:domain-decomp}. In particular, we attain iteration counts for the Signorini problem that are robust to increasing incompressibility of the elastic material, i.e.~as the Lam\'e parameter $\lambda \to \infty$.

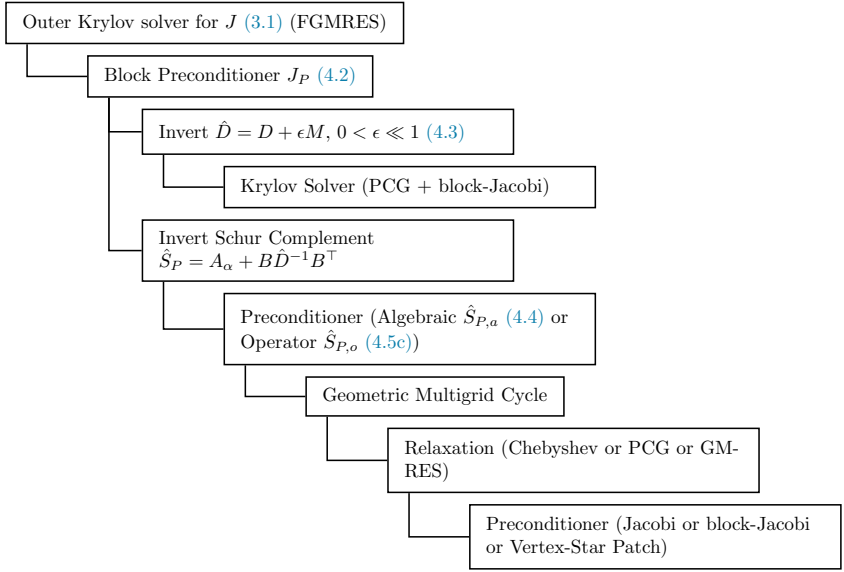
\begin{figure}[ht]
\centering
\begin{tikzpicture}[
    box/.style={draw=black, thick, align=left,
        inner xsep=3mm, inner ysep=2mm},
    wrap/.style={box, text width=6.2cm},
    connect/.style={thick},
]
\node[box] (outer) {Outer Krylov solver for $J$ \cref{eq:saddle-point} (FGMRES)};
\node[box, anchor=north west]
    (block) at ([xshift=15mm,yshift=-2mm]outer.south west)
    {Block Preconditioner $J_P$ \cref{eq:schur2}};
\node[wrap, anchor=north west]
    (d) at ([xshift=10mm,yshift=-2mm]block.south west)
    {Invert $\hat{D}=D+\epsilon M$, $0<\epsilon\ll 1$ \cref{eq:D-hat}};
\node[wrap, anchor=north west]
    (krylov) at ([xshift=15mm,yshift=-2mm]d.south west)
    {Krylov Solver (PCG $+$ block-Jacobi)};
\node[wrap, anchor=north west]
    (schur) at
    ($(block.south west |- krylov.south west)+(10mm,-2mm)$)
    {Invert Schur Complement \\
     $\hat{S}_P=A_\alpha+B\hat{D}^{-1}B^\top$};
\node[wrap, anchor=north west]
    (schurprec) at ([xshift=15mm,yshift=-2mm]schur.south west)
    {Preconditioner (Algebraic $\hat{S}_{P,a}$ \cref{eq:Sd} or Operator $\hat{S}_{P,o}$ \cref{eq:So})};
\node[box, anchor=north west]
    (mg) at ([xshift=15mm,yshift=-2mm]schurprec.south west)
    {Geometric Multigrid Cycle};
\node[wrap, anchor=north west]
    (relax) at ([xshift=15mm,yshift=-2mm]mg.south west)
    {Relaxation (Chebyshev or PCG or GMRES)};
\node[wrap, anchor=north west]
    (patch) at ([xshift=15mm,yshift=-2mm]relax.south west)
    {Preconditioner (Jacobi or block-Jacobi or Vertex-Star Patch)};

\draw[connect] ([xshift=4mm]outer.south west) |- (block.west);
\draw[connect] ([xshift=4mm]block.south west) |- (d.west);
\draw[connect] ([xshift=4mm]d.south west) |- (krylov.west);
\draw[connect] ([xshift=4mm]block.south west) |- (schur.west);
\draw[connect] ([xshift=4mm]schur.south west) |- (schurprec.west);
\draw[connect] ([xshift=4mm]schurprec.south west) |- (mg.west);
\draw[connect] ([xshift=4mm]mg.south west) |- (relax.west);
\draw[connect] ([xshift=4mm]relax.south west) |- (patch.west);
\end{tikzpicture}
\caption{Block preconditioning strategy for solving the PG Newton systems \cref{eq:saddle-point}.}\label{fig:block}
\end{figure}

We summarize the block preconditioning strategy in \cref{fig:block}.

\subsection{Subspace correction methods}
\label{sec:domain-decomp}
In this section, we introduce the vertex-star patch and Vanka-type iteration as \emph{additive Schwarz methods} (ASM) \cite{pcpatch}.

Let $Z_h$ denote a (possibly mixed) finite element space on the domain $\Omega$ and consider the discretized variational problem, find $z_h \in Z_h$ satisfying
\begin{align}
\label{eq:asm1}
a(z_h, w_h) = F(w_h) \quad \forall w_h \in Z_h
\end{align}
where $a : Z_h \times Z_h \to \mathbb{R}$ and $F : Z_h \to \mathbb{R}$ denote a bilinear and a linear form, respectively. We attempt to solve \cref{eq:asm1} by solving several smaller local variational problems and combining the local corrections. In this vein, consider the decomposition of the FEM space $Z_h = \sum_i Z_h^{(i)}$, where the sum need not be a direct sum. Fix an iterate $z_{j,h}$, $j \in \mathbb{N}$ approximating $z_h$ and consider the corresponding local subspace correction variational problem: find $e_h^{(i)} \in Z_h^{(i)}$ satisfying
\begin{align}
\label{eq:asm2}
a(e^{(i)}_h, w_h) = F(w_h) - a(z_{j,h},w_h) \quad \forall w_h \in Z^{(i)}_h.
\end{align}
For each $i$, we extend $e^{(i)}_h \in  Z^{(i)}_h$ by zero so it is well-defined as a function $e^{(i)}_h \in Z_h$ on all of $\Omega$. Then the ASM iteration is $z_{j+1,h} = z_{j,h} + \sum_i e^{(i)}_h$. Note that the local problems \cref{eq:asm2} are independent for each $i$ and can be solved in parallel. The vertex-star iteration and Vanka-type iteration are both an ASM with a specific choice for $Z_h$ and decomposition $\{Z^{(i)}_h\}$.

\subsubsection{Vertex-star patch iteration}
The vertex-star patch iteration is an ASM only applied to variational equations posed over the primal FEM space $V_h$. Thus $Z_h = V_h$. The decomposition is based on the mesh entities. For each vertex $v_i$ in the mesh, each local subspace is
\[
Z^{(i)}_h = \{ \varphi_j \in V_h : v_i \in \operatorname{int}(\supp(\varphi_j)) \},
\]
where $\mathrm{int}$ denotes the strict interior of a set and $\supp$ the support of a function. In other words, $Z^{(i)}_h$ consists of the FEM basis functions for $V_h$ that include the vertex $v_i$ within the interior of their support. The collection of DoFs around each vertex is known as a patch.  In \cref{fig:patch:mtw,fig:patch:cg2,fig:patch:cg2a}, we demonstrate the patches for an $\mathrm{MTW}$, $\mathrm{CG}_2$ and $\mathrm{CG}_2$-A discretization. Note that the vertex-star patch iteration for a $\mathrm{CG}_1$ element is equivalent to the Jacobi iteration for a scalar-valued $V_h$ (or block-Jacobi for a vector-valued $V_h$). For all the problems considered in this work, the resulting local variational problems \cref{eq:asm2} are SPD.

Vertex-star patch iterations are a key ingredient in designing parameter-robust solvers for singularly perturbed problems. In particular, using vertex-star patch iterations as the relaxation method in a multigrid cycle can lead to an iterative solver with robust iteration counts for nearly-incompressible linear elasticity \cite{schoberl1999,farrell2022}. We shall observe in \Cref{sec:examples} that the same mechanism leads to robust iteration counts for nearly-incompressible Signorini. The second key ingredient is a specialized prolongation operator whose continuity constant is independent of the singular perturbation parameter (in this case the Lam\'e parameter $\lambda$). For the MTW element, the standard prolongation operator suffices whereas for the CG-A element a non-standard prolongation is required, cf.~\cite[Sec.~3]{farrell2022}.

\subsubsection{Vanka-type iteration}
\begin{figure}[h!]
\centering
\subfloat[MTW (Star) \label{fig:patch:mtw}]{\includegraphics[width =0.24\textwidth]{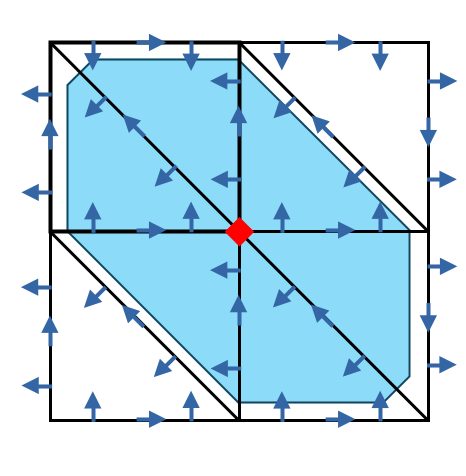}}
\subfloat[$\mathrm{CG}_2$ (Star) \label{fig:patch:cg2}]{\includegraphics[width =0.24\textwidth]{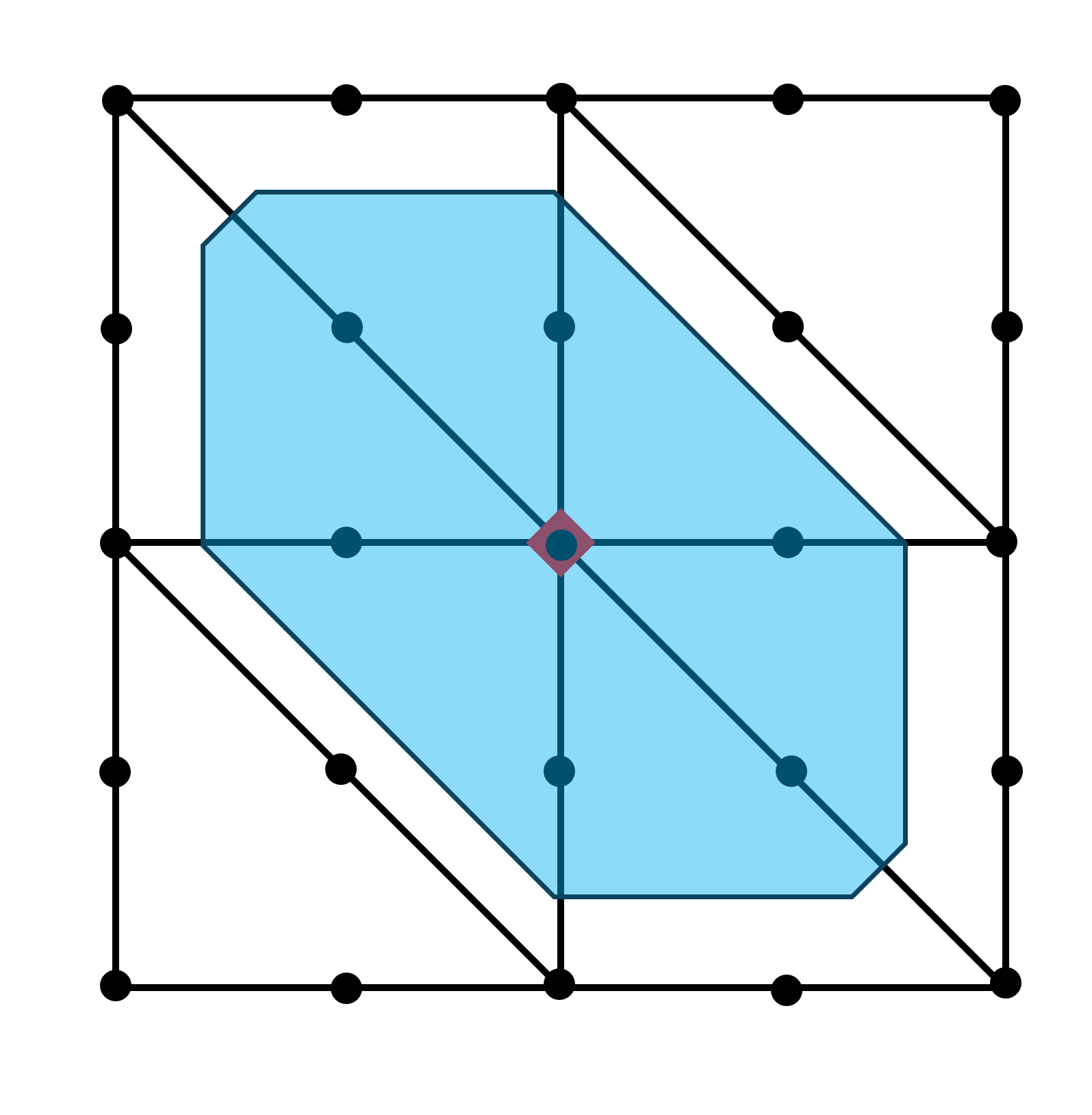}}
\subfloat[$\mathrm{CG}_2$-A (Star) \label{fig:patch:cg2a}]{\includegraphics[width =0.24\textwidth]{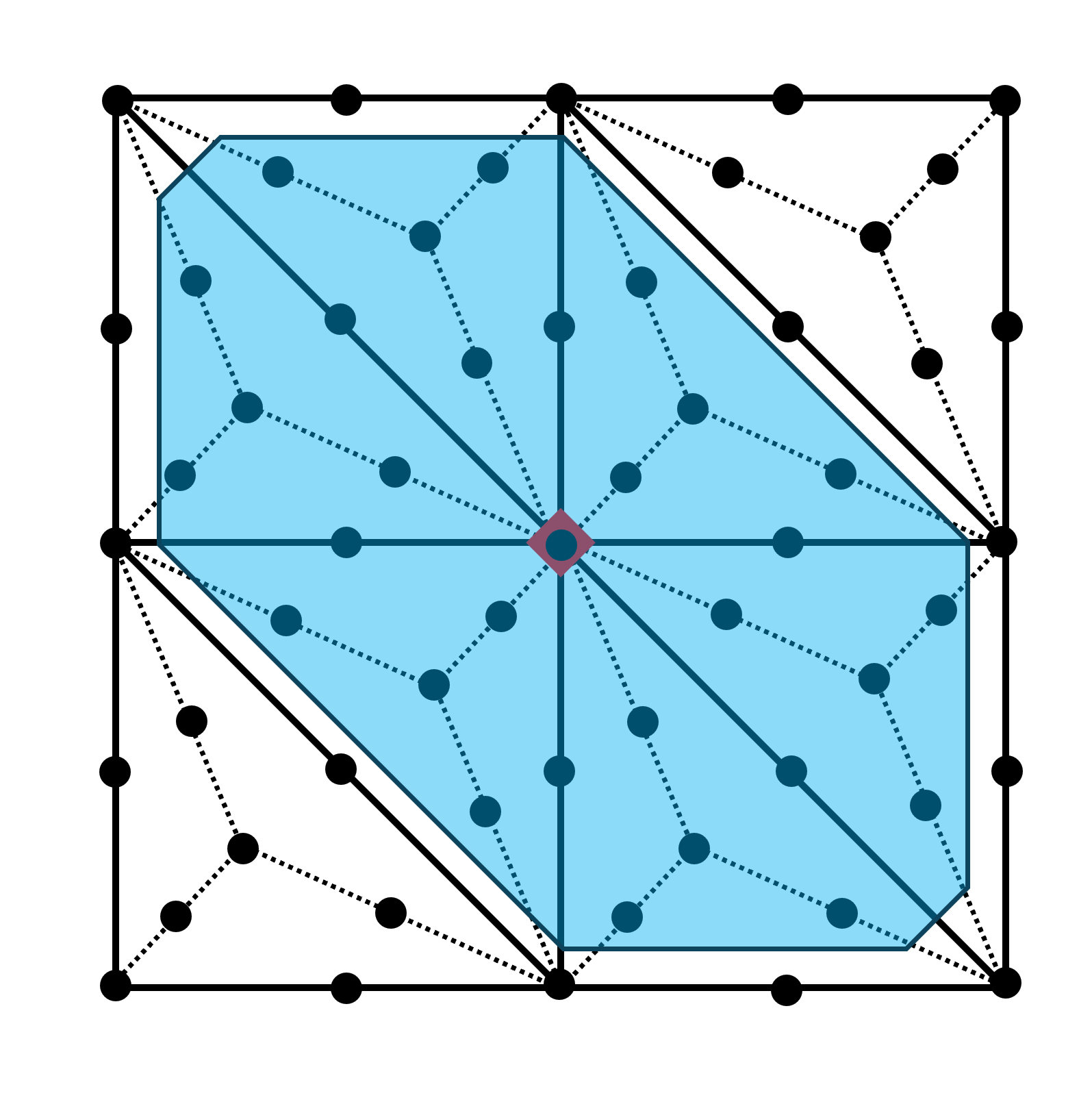}}
\subfloat[$\mathrm{CG}_1\times\mathrm{CG}_1$ (Vanka)  \label{fig:patch:vanka}]{\includegraphics[width =0.24\textwidth]{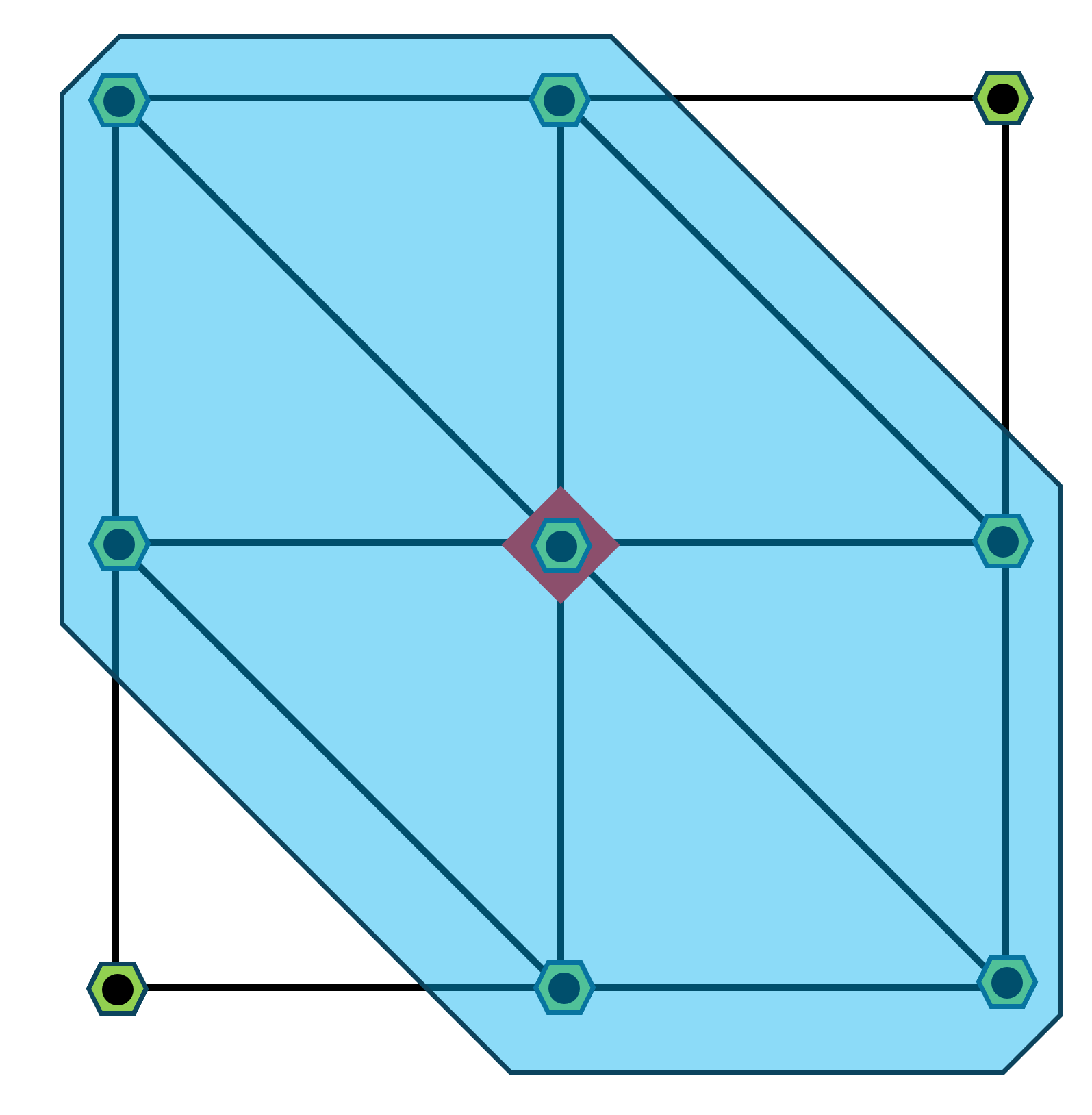}}
\caption{(a)-(c) depict the vertex-star patch around an internal vertex (red diamond) on a two-dimensional mesh for the MTW (18 DoFs), $\mathrm{CG}_2$ (scalar valued: 7 DoFs, vector-valued: 14 DoFs) and $\mathrm{CG}_2$-A (scalar-valued: 31 DoFs, vector-valued: 62 DoFs) elements. The DoFs are depicted by black dots or blue arrows. In (d), the Vanka-type patch is depicted around an internal vertex for a $\mathrm{CG}_1 \times \mathrm{CG}_1$ discretization. The black dots and green hexagons indicate the DoFs for each field, respectively. Each patch contains 14 DoFs.}
\label{fig:patch}
\end{figure}

The Vanka-type iteration is an ASM applied to variational equations posed over the mixed function space for $(u,\psi)$, i.e.~$Z_h = V_h \times Q_h$. Again the decomposition is vertex-based and $Z^{(i)}_h = V^{(i)}_h \times Q^{(i)}_h$ where
\begin{align*}
V^{(i)}_h = \{ \varphi_j \in V_h : v_i \in \operatorname{cl}(\supp(\varphi_j)) \} \;\;\text{and} \;\; Q^{(i)}_h = \{ \eta_j \in Q_h : v_i \in \operatorname{cl}(\supp(\eta_j)) \},
\end{align*}
where $\mathrm{cl}$ denotes the closure of a set. Note that Vanka-type patches include DoFs in the closure of the support of the basis functions rather than just the interior. The Vanka-type patches considered in this work are non-standard in the sense that typically one does not include the closure of the supports for the ``Lagrange-multiplier" space $Q_h$. 
This choice is due to the different saddle point structure of our problems.
In \cref{fig:patch:vanka}, we demonstrate the Vanka-type patch for a $\mathrm{CG}_1 \times \mathrm{CG}_1$ discretization. Unlike the vertex-star patch iteration, the local variational problems \cref{eq:asm2} are saddle point problems.

\subsection{Monolithic multigrid}
\label{sec:monolithic}
In this subsection, we briefly introduce another preconditioner candidate for the PG Newton systems. We compare this preconditioner with the block preconditioner of \Cref{sec:block} in \Cref{sec:examples}.

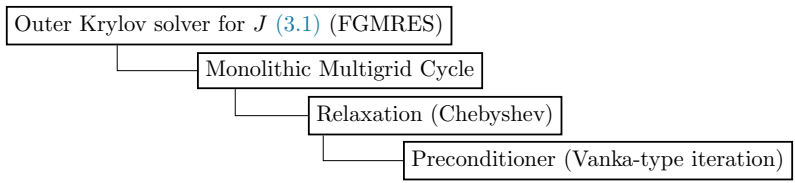
\begin{figure}[ht]
\centering
 \begin{tikzpicture}[%
every node/.style={draw=black, thick, anchor=west},
grow via three points={one child at (-0.5,-0.7) and
    two children at (-0.5,-0.7) and (-0.5,-1.4)},
edge from parent path={(\tikzparentnode.190) |- (\tikzchildnode.west)}]
\node {Outer Krylov solver for $J$ \cref{eq:saddle-point} (FGMRES)}
child { node {Monolithic Multigrid Cycle}
    child { node {Relaxation (Chebyshev)}
     			child { node{Preconditioner (Vanka-type iteration)}}}
};
\end{tikzpicture}
\caption{Monolithic multigrid strategy for solving the PG Newton systems \cref{eq:saddle-point}.}\label{fig:vanka}
\end{figure}

Rather than construct a block preconditioner for the FGMRES method, one can directly approximate the inverse of the Jacobian in \cref{eq:saddle-point} via monolithic geometric multigrid applied directly to the full Jacobian. However, due to the saddle point structure, traditional relaxation methods, such as Chebyshev preconditioned by a Jacobi iteration, are ineffective and the multigrid cycle does not converge. Thus we borrow a technique from other saddle point literature and instead use the Vanka-type iteration as the relaxation \cite{farrell2021b,vanka1986,voronin2025}. More specifically, we use it as a preconditioner for a handful of Chebyshev iterations on each level of the cycle.  The monolithic multigrid strategy for solving \cref{eq:saddle-point} is summarized in \cref{fig:vanka}.

In \Cref{sec:examples}, we see that monolithic multigrid is an effective strategy for the obstacle problem. However, in unreported experiments, monolithic multigrid preconditioning for PG would diverge when applied to the strain-constrained elasticity problem. Moreover, an implementation of a Vanka-type iteration for the Signorini problem is nontrivial since the latent variable is only supported on the subdomain on the boundary.

\section{Numerical experiments}
\label{sec:examples}

In this section, we study an instance of an obstacle, Signorini, and strain-constrained problem in both two and three dimensions. All 2D simulations are run in serial and the 3D simulations are distributed across 32 cores on a workstation with up to 2048GB of RAM. The tables are deferred to the end of the section for readability.

\subsection*{Data availability}  Software to generate all the tables and figures in this article may be found at \cite{github:proximalgalerkin}. The specific version used in this paper is archived on Zenodo \cite{zenodo}. The implementation is written in the open-source finite element package Firedrake \cite{firedrake} using PETSc \cite{petsc} as the linear algebra backend. The vertex-star patch and Vanka-type iterations are implemented using PCPatch  \cite{pcpatch}. The meshes are either natively constructed in Firedrake or via Netgen \cite{netgen}. \Cref{sec:examples:signorini} utilizes the recently introduced \texttt{Submesh} functionality in Firedrake to construct the latent variable finite element space on the boundary \cite{sagiyama2025}. The algebraic Schur approximant is constructed via the \texttt{Slate} functionality \cite{slate}. Visualizations are made in ParaView \cite{paraview}.

\subsection*{Notation}  $N$ and $p$ parametrize the maximum mesh size $h=1/N$ and degree of the FEM, respectively. The ``levels" column denotes the number of levels in a multigrid cycle including the base mesh. The ``PG" and ``Newton" columns record the number of proximal steps and average Newton iterations per proximal step to reach convergence, respectively. The FGMRES column denotes the average number of FGMRES iterations per Newton step. Each FGMRES iteration requires a single multigrid cycle. In the same FGMRES column, we include, in brackets, the maximum number of FGMRES iterations a Newton step required over the full PG solve. In the case of a block preconditioner, each FGMRES iteration requires two applications of the inverse of the preconditioner's latent block $-\hat{D}$. The ``PCG" column denotes the average number of block-Jacobi PCG iterations to invert $\hat{D}$ per action of $\hat{D}^{-1}$. Again, the bracketed numbers are the maximum number of PCG iterations a $\hat{D}$-inverse took over the whole PG run.

For the linear solvers, we use ``LU" to denote a sparse LU factorization of the full Newton system. ``MG-Vanka" will denote using the monolithic multigrid strategy outlined in \cref{fig:vanka}. For the block preconditioners, we use ``AS-MG" and ``OS-MG" to denote using the algebraic Schur complement approximant $\hat{S}_{P,a}$ \cref{eq:Sd} and operator Schur complement approximant $\hat{S}_{P,o}$ \cref{eq:So}, respectively. We then append these names with ``-Ch'', ``-PCG'', ``-GMRES'' to distinguish a Chebyshev, conjugate gradient or GMRES relaxation iteration followed by ``-J", ``-bJ" or ``-S" for a Jacobi, block-Jacobi or vertex-star patch preconditioner of the relaxation iteration. We summarize the linear solvers used in this section  in \cref{tab:solvers}.

\begin{table}[h!]
\renewcommand{\arraystretch}{1.2}
\centering
\scalebox{0.9}{
\begin{tabular}{|l|l|}
\hhline{|-|-|}
\rowcolor{lightgray!10}
\multicolumn{1}{c|}{Solver} & \multicolumn{1}{c|}{Description} \\
\hline
LU & LU factorization of the Newton matrix $J$ in \cref{eq:saddle-point}. \label{row:lu}\\
\hline
MG-Vanka & The monolithic multigrid strategy outlined in \cref{fig:vanka}. \label{row:vanka}\\
\hline
\multirow{4}{*}{AS-MG-Ch-J}
 & Block preconditioner strategy of \cref{fig:block}. The  Schur \\
& complement is approximated by inverting the algebraic Schur \\
& complement $\hat{S}_{P,a}$ in \cref{eq:Sd} with geometric multigrid. The   \\
& relaxation is Chebyshev preconditioned by a Jacobi iteration. \\ 
\hline
\multirow{2}{*}{AS-MG-Ch-S}
 & Same as AS-MG-Ch-J but the relaxation is Chebyshev  \\
 &  preconditioned   by a vertex-star patch iteration. \\ 
 \hline
 \multirow{2}{*}{AS-MG-GMRES-bJ}
 & Same as AS-MG-Ch-J but the relaxation is GMRES    \\
 &  preconditioned  by a block-Jacobi iteration. \\ 
 \hline
\multirow{4}{*}{OS-MG-Ch-J}
 & Block preconditioner strategy of \cref{fig:block}. The  Schur \\
& complement is approximated by inverting the operator Schur \\
& complement $\hat{S}_{P,o}$ in \cref{eq:So} with geometric multigrid. The \\ 
& relaxation  is Chebyshev  preconditioned by a Jacobi iteration. \\ 
\hline
\multirow{2}{*}{OS-MG-Ch-S}
 & Same as OS-MG-Ch-J but the relaxation is Chebyshev   \\
 &  preconditioned   by a vertex-star patch iteration. \\ 
\hline
\multirow{2}{*}{OS-MG-PCG-S}
 &Same as OS-MG-Ch-J but the relaxation is PCG     \\
 & preconditioned by a vertex-star patch iteration. \\ \hline
  \multirow{2}{*}{OS-MG-GMRES-bJ}
 & Same as OS-MG-Ch-J but the relaxation is GMRES    \\
 &   preconditioned by a block-Jacobi iteration. \\ 
\hline
\end{tabular}}
\caption{A summary of the linear solvers for the Newton linear systems \cref{eq:saddle-point} investigated in \Cref{sec:examples}.}
\label{tab:solvers}
\end{table}

\subsection{Obstacle problem}
\label{sec:examples:obstacle}
In this example, we consider an obstacle problem \cref{eq:obstacle} posed on an L-shaped domain in 2D and a box with a corner removed in 3D. More specifically:
\begin{align}
\Omega = (0,2)^d \backslash (1,2)^d, \quad f \equiv 20, \quad \phi \equiv 1.
\end{align}
In both dimensions, we initialize the base mesh and refine cells  within a radius of 0.05 around the re-entrant corner. Uniform mesh refinements for the geometric multigrid hierarchy are subsequently built upon this base mesh. We discretize the problem with equal-order continuous Lagrange finite elements for $u$ and $\psi$ ($\mathrm{CG}_p \times \mathrm{CG}_p$, $p \in \{1,2,3\}$).

\begin{figure}[h!]
\centering
\includegraphics[width =0.28 \textwidth]{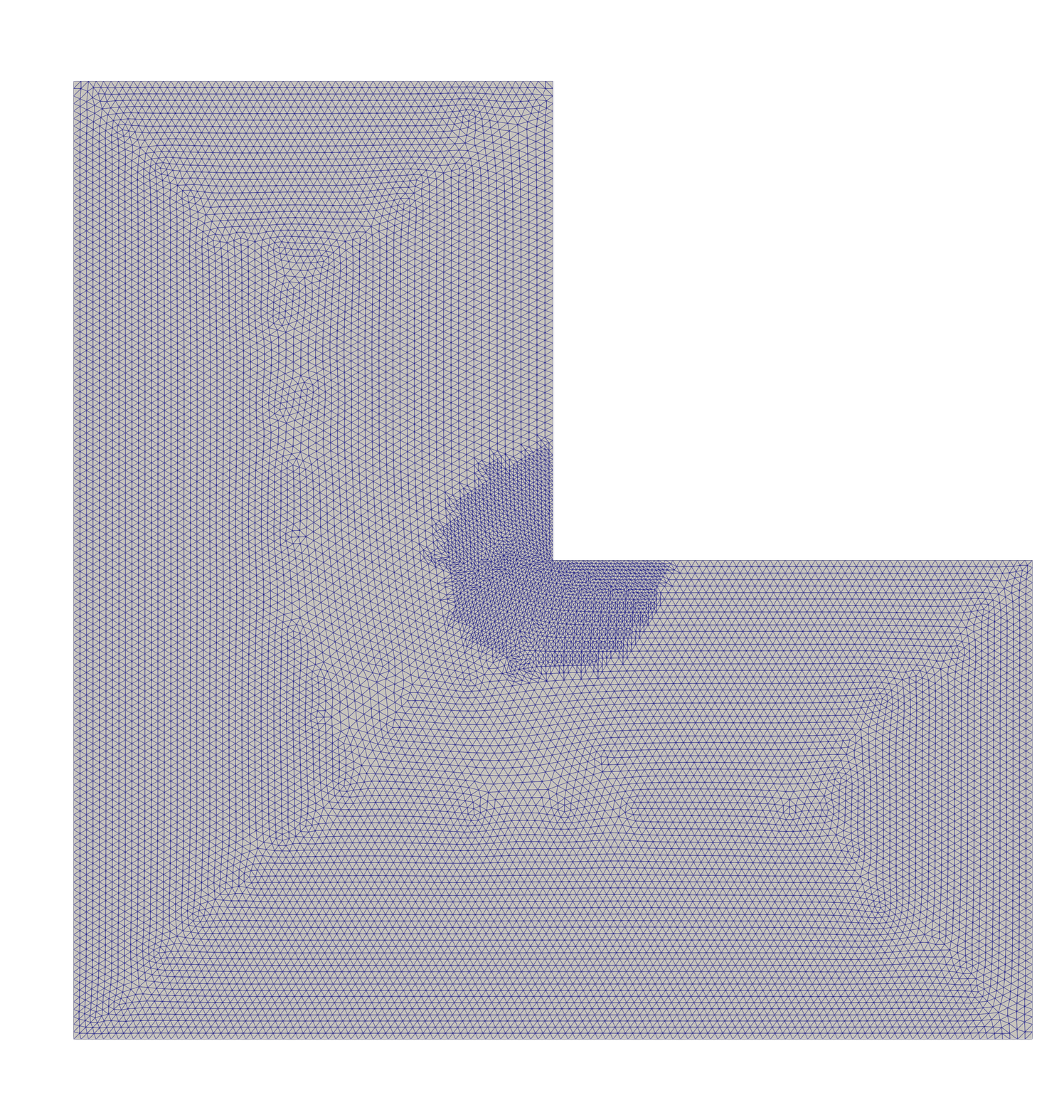} 
\includegraphics[width =0.36 \textwidth]{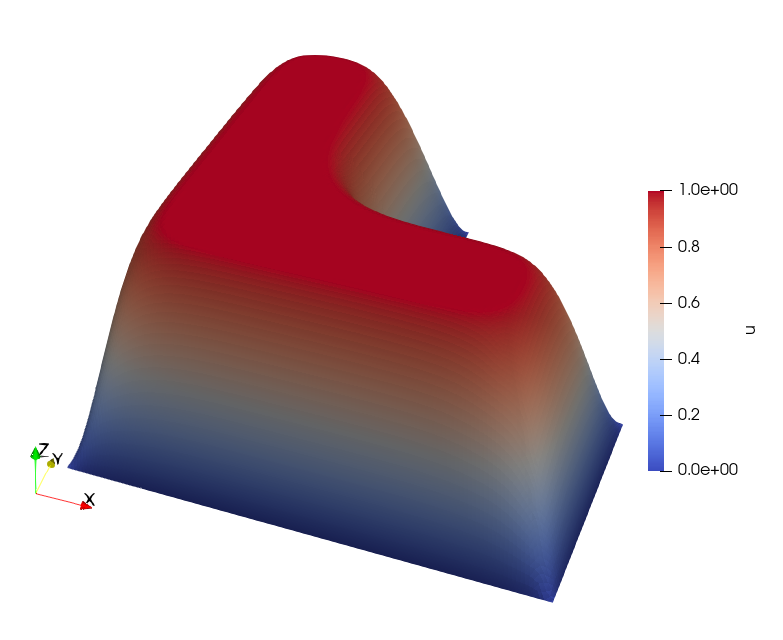}
\includegraphics[width =0.33 \textwidth]{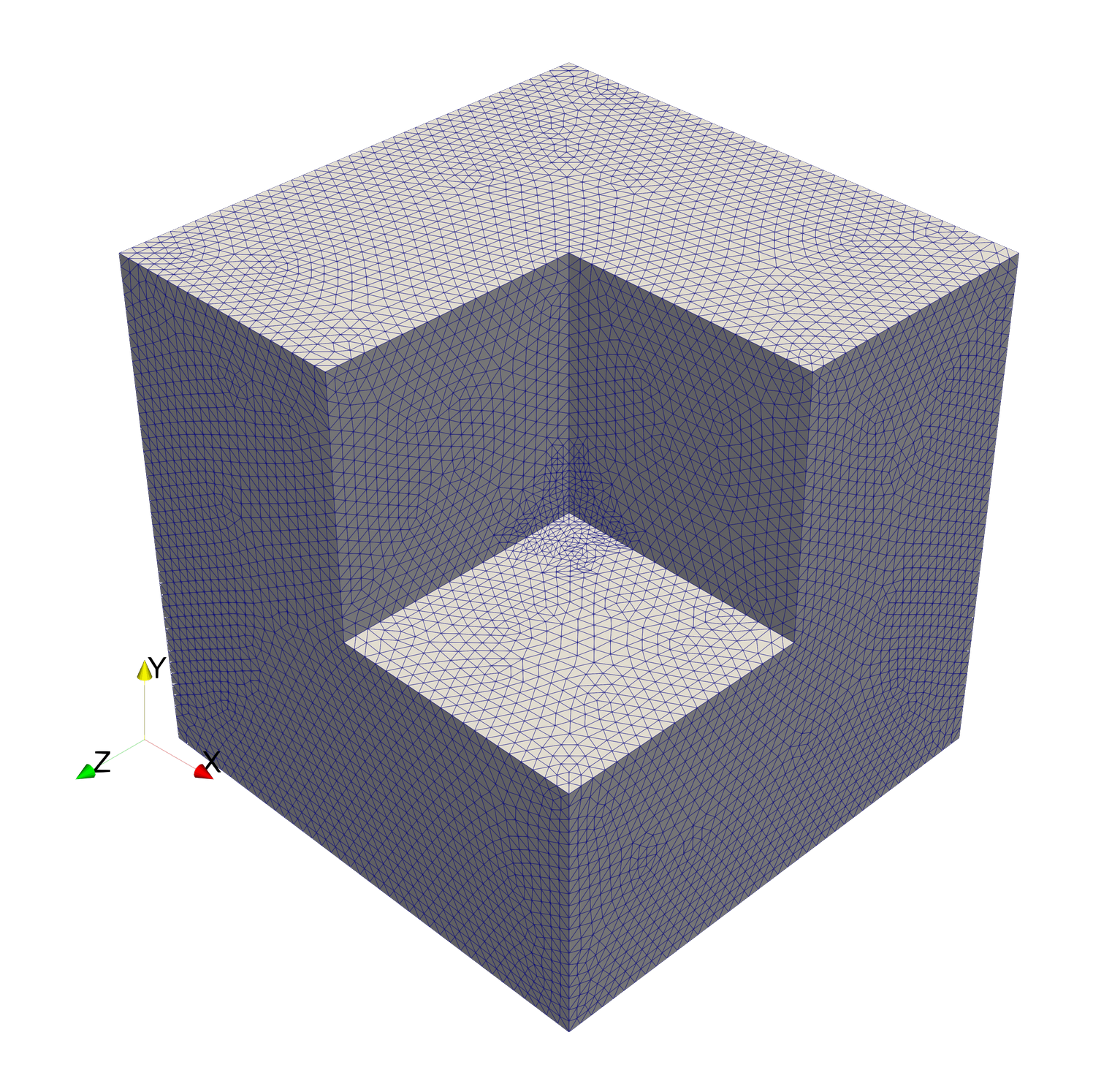}
\caption{(Left \& Right) Mesh of the 2D ($N=64$) and 3D ($N=20$) domains in the obstacle problem of \Cref{sec:examples:obstacle}. (Middle) Solution of the 2D obstacle problem.}
\label{fig:obstacle-mesh}
\end{figure}

All runs are initialized with $\alpha_0 = 10^{-1}$, use the proximal update rule $\alpha_k = \sqrt{2} \alpha_{k-1}$, and terminate once $\|u^k-u^{k-1}\|_{H^1(\Omega)} \leq 10^{-3}$. We pick $u^0 = \psi^0 = 0$. The initial iterate for Newton's method at the $k^\text{th}$-proximal step is $(u^{k-1}, \psi^{k-1})$. The Newton solves terminate when the $\ell^2$-norm of the residual drops below $10^{-5}$ and we use an absolute and relative tolerance of $10^{-6}$ for the Krylov solvers. We apply five relaxation iterations on each level of the multigrid cycle, except on the coarsest level where an LU factorization is used.

In \cref{fig:obstacle-mesh} we visualize a choice of base mesh and re-entrant corner refinement in 2D (left) and 3D (right). We also plot the 2D solution (middle).

\subsubsection{2D Obstacle}
In \cref{tab:obstacle2d-1}, we investigate a sparse LU factorization to solve every Newton system together with a monolithic multigrid preconditioner with a Vanka-type relaxation (MG-Vanka, \cref{fig:vanka}). We report the findings on three different meshes, obtained from uniformly refining a base mesh once, twice, and thrice, as well as three choices of polynomial degree ($p= 1$, 2 and 3). For these linear solvers, a $D$-inverse is not needed. Hence we set the preconditioner perturbation to zero, $\epsilon=0$. For MG-Vanka, we examine two- and three-level multigrid cycles.

All the PG solves converged in 9 proximal steps independent of the discretization. The average number of Newton iterations per proximal step was under 1.8, dropping to 1.44 on the finest grid with $p=3$. On average, fewer than three FGMRES iterations were required per Newton step and, therefore, under three multigrid cycles. These average iteration counts also decreased as the discretization was refined.

In \cref{tab:obstacle2d-2}, we report the effectiveness of four different realizations of the block preconditioning strategy in \cref{fig:block}: OS-MG-Ch-J, OS-MG-Ch-S, AS-MG-Ch-J, and AS-MG-Ch-S.
We fix the perturbation of the preconditioner's $\hat{D}$-block to $\epsilon = 10^{-5}$ and test two- and three-level cycles on two different meshes and discretization degree $p \in \{1,2,3\}$.

With the block preconditioner, all the PG solves converged in 9 proximal steps. The average number of Newton solves per proximal step was identical across the four solvers peaking at 1.78 and improving to 1.44 as the discretization was refined. Concerning the average number of outer FGMRES iterations per Newton step, the OS-MG-Ch-S solver performed overall best, followed by OS-MG-Ch-J, AS-MG-Ch-S, and finally AS-MG-Ch-J. The best average was 1.77 FGMRES per Newton step and the worst was 4.94. The iteration counts \emph{improved} as the discretization was refined. The average number of PCG iterations per $\hat{D}$-inverse remained bounded with a best average iteration count of 3.95 and a worst of 5.38.

Note that although the OS-MG-Ch-S solver produced the smallest average FGMRES iteration counts, a Jacobi iteration is cheaper than a vertex-star patch iteration. Hence, we recommend the OS-MG-Ch-J solver for performance.

\subsubsection{3D Obstacle}
 
In \cref{tab:obstacle3d}, we compare the monolithic multigrid strategy MG-Vanka against OS-MG-Ch-J, OS-MG-Ch-S, AS-MG-Ch-J, and AS-MG-Ch-S. We fix a base mesh and consider two-level and three-level multigrid cycles built on uniform refinements of this mesh. We also report the iteration counts for the discretization degrees $p=1$ and $p=2$. We set $\epsilon=0$ in the MG-Vanka strategy, 
and $\epsilon=10^{-4}$, otherwise.

All PG solves converged in 25 proximal steps with an average of 1.6 Newton iterations per proximal step or less. MG-Vanka produced the lowest average FGMRES iterations per Newton step with a high of 2.62 and low of 2.18. The next best was OS-MG-Ch-S followed by OS-MG-Ch-J, AS-MG-Ch-S and finally AS-MG-Ch-J. Note that out of the block preconditioning strategies, only OS-MG-Ch-J was able to run on our machine for the finest discretization. The highest average FGMRES iteration count was 4.05 and lowest was 2.53. The average PCG iteration count per $\hat{D}$-inverse was under 8 when $p=1$ and 15 when $p=2$.

As in 2D, although patch relaxation produced the best iteration counts, the Jacobi iteration is significantly cheaper. Hence, we recommend the OS-MG-Ch-J solver for the best performance.

\subsection{Signorini problem}
\label{sec:examples:signorini}
In this subsection we fix an example of the Signorini problem \cref{eq:signorini}. The goal is to not only examine the robustness of the preconditioned PG strategy to the mesh and degree but also its robustness as $\lambda \to \infty$, i.e.~with respect to the incompressibility of the elastic material. We shall observe that, via vertex-star patch relaxation, one can attain robustness across all levels: the number of proximal steps, Newton steps, and Krylov iterations to solve each Newton system.

\begin{figure}[h!]
\centering
\subfloat[2D: $\lambda=10^5$ ($\nu=0.49998$) \label{fig:signorini:2d}]{\includegraphics[width = 0.9\textwidth]{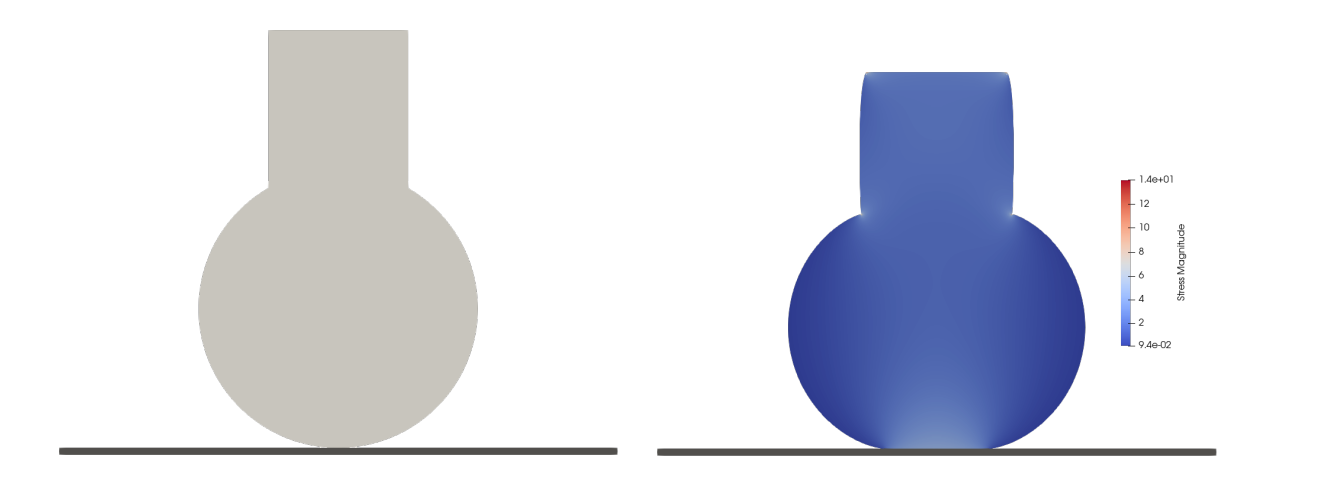}}\\
\subfloat[3D: $\lambda=10^4$ ($\nu=0.49$) \label{fig:signorini:3d}]{\includegraphics[width =0.9\textwidth]{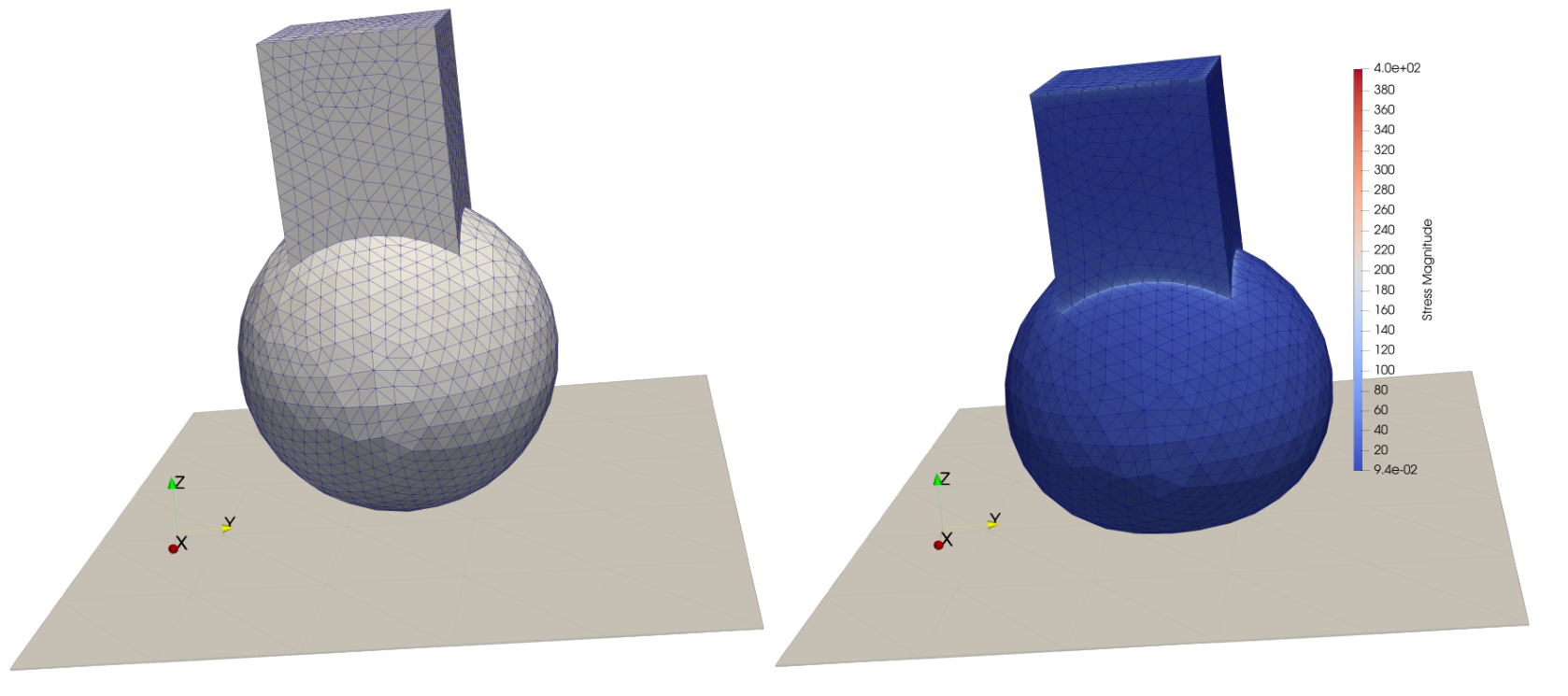} }
\caption{The undeformed and deformed elastic bodies described in \Cref{sec:examples:signorini}. Here $\nu$ denotes the Poisson's ratio and the discretization is an Alfeld-split $(\mathrm{CG}^A_d)^d \times \mathrm{CG}_d$ FEM.}
\label{fig:signorini}
\end{figure}

In 2D we choose the model parameters:

\begin{align}
\begin{gathered}
\Omega = \left\{ (x_1,x_2) \in \mathbb{R}^2 : x_1^2+x_2^2< 1 \right\} \cup [(-0.5,0.5) \times (0.5,2)], \\
f = (0,0)^\top, \quad
\mu = 5, \quad
\Gamma_C = \{x_2 < 0\}, \quad n = (0,-1)^\top,\\
 \phi = 1-\sqrt{1-x_1^2}, \quad \Gamma_D = \{x_2=2\}, \quad g=(0,-0.3)^\top.
\end{gathered}
\end{align}
The Lam\'e parameter $\lambda$ will be varied in the experiments. This setup models a keyhole domain where the bottom of the disk section rests on a surface. The top surface of the domain is compressed vertically downwards by 0.3 and the problem is to find the resulting deformation of the elastic body as it compresses against the obstacle surface. The solution is visualized in \cref{fig:signorini:2d} when $\lambda = 10^5$. The relationship with the Poisson's ratio is $\nu = \frac{\lambda}{2\lambda + 2\mu}$.

The three-dimensional example has the following setup:
\begin{align}
\begin{gathered}
\Omega = \left\{ (x_1,x_2,x_3) \in \mathbb{R}^3 : x_1^2+x_2^2 + x_3^2< 1 \right\} \cup [(-0.5,0.5)^2 \times (0.5,2)], \\
f = (0,0,0)^\top, \quad
\mu = 200, \quad
\Gamma_C = \{x_3 < 0\}, \quad n = (0,0,-1)^\top,\\
\quad \phi = 1-\sqrt{1-x_1^2-x_2^2}, \quad \Gamma_D = \{x_3 = 2\}, \quad g = (0,0,-0.3)^\top.
\end{gathered}
\end{align}
The solution when $\lambda = 10^4$ is visualized in \cref{fig:signorini:3d}.

\subsection{2D Signorini}

All runs in 2D are initialized with $\alpha_0 = 1$, use the proximal update rule $\alpha_k = \min\{2 \alpha_{k-1}, 10^2\}$, and terminate once $\|u^k-u^{k-1}\|_{H^1(\Omega)} \leq 10^{-3}$. We pick $u^0 = \psi^0 = 0$ where the Dirichlet boundary conditions are pre-applied to $u^0$ before the first Newton solve. The initial iterate for Newton's method at the $k^\text{th}$-proximal step is $(u^{k-1}, \psi^{k-1})$. The Newton solves terminate when the $\ell^2$-norm of the residual drops below $10^{-5}$ and we use an absolute and relative tolerance of $10^{-6}$ for the Krylov solvers. We apply five relaxation iterations on each level of the MG schemes except on the coarsest level where an LU factorization is used.  

In \cref{tab:signorini2d-1} we fix $\lambda=10$ and investigate the effectiveness of an LU factorization, and the block preconditioner \cref{fig:block} with the operator Schur complement approximant \cref{eq:So} and Chebyshev relaxation with a Jacobi preconditioner (OS-MG-Ch-J) or a vertex-star patch preconditioner (OS-MG-Ch-S). We consider three discretizations of $(u,\psi)$: equal-order continuous FEM $(\mathrm{CG}_p)^2 \times \mathrm{CG}_p$, $p \in \{1,2\}$, denoted CG, an Alfeld-split quadratic FEM for $u$ $(\mathrm{CG}^A_2)^2 \times \mathrm{CG}_2$, denoted CG-A, and the non-conforming FEM $\mathrm{MTW}_1 \times \mathrm{DG}_0$, denoted MTW. We fix the perturbation of the block preconditioner's $\hat{D}$-block as $\epsilon =10^{-4}$. For the CG and MTW discretizations, we use the default prolongation operator, whereas a specialized prolongation is used for CG-A \cite[Sec.~3]{farrell2022}.

All the PG solves converged in six proximal steps averaging between 2.50 and 2.67 Newton iterations per proximal step. 

We applied a two- and three-level OS-MG-Ch-J solver to the  $(\mathrm{CG}_p)^2 \times \mathrm{CG}_p$, $p \in \{1,2\}$, discretization. We observe robust FGMRES iteration counts with respect to the mesh size, number of MG levels and discretization degree $p$, averaging between 4.19 and 4.62 FGMRES iterations per Newton step with a maximum number of 9 FGMRES iterations. The number of PCG iterations to invert the $\hat{D}$-block also remained robust to the mesh size and degree, averaging between 1.55 and 1.70 PCG iterations per inverse.

We then investigated a two- and three-level OS-MG-Ch-S solver for all three discretizations. Once again, the FGMRES iteration counts were mostly robust to mesh size and discretization degree. The lowest iteration counts were achieved via the CG-A discretization, whereas the MTW discretization produced the highest iteration counts. Again the PCG iterations remain small. The MTW discretization inverted the $\hat{D}$-block with one PCG iteration per inverse since $\hat{D}$ is diagonal when using a $\mathrm{DG}_0$ discretization for $\psi$.

In \cref{tab:signorini2d-2}, we test the robustness of the two-level OS-MG-Ch-S solver as $\lambda \to \infty$. We fix a mesh for each discretization and solely consider $p=1$ for the CG discretization. We pick the perturbation $\epsilon=10^{-4}$. For the CG discretization, the PCG iterations for the $\hat{D}$-inverse remain similar, but the FGMRES iteration counts increase as $\lambda$ increases averaging 3.88, 8.53, and 26.33 when $\lambda = 10, 10^2$ and $10^3$, respectively. When $\lambda=10^4$ the solver fails to converge. In contrast, the iteration counts for the CG-A and MTW discretizations remain robust. For CG-A, the average FGMRES iterations per Newton step increases from 3.31 when $\lambda=10$ to 4.92 when $\lambda=10^5$. Similarly for the MTW discretization it increases from 9.69 to 10.00.

We conclude that the OS-MG-Ch-S solver results in iteration counts robust to the incompressibility of the material when discretized with a CG-A or MTW discretization.

\subsection{3D Signorini}

All runs in 3D are initialized with $\alpha_0 = 10^{-1}$, use the proximal update rule $\alpha_k = \min\{\sqrt{2} \alpha_{k-1}, 10\}$, and terminate once $\|u^k-u^{k-1}\|_{H^1(\Omega)} \leq 10^{-3}$. We pick $u^0 = \psi^0 = 0$  where the Dirichlet boundary conditions are pre-applied to $u^0$ before the first Newton solve. The initial iterate for Newton's method at the $k^\text{th}$-proximal step is $(u^{k-1}, \psi^{k-1})$. The Newton solves terminate when the $\ell^2$-norm of the residual drops below $10^{-5}$ and we use an absolute and relative tolerance of $10^{-6}$ for the Krylov solvers. We apply five relaxation iterations on each level of the MG schemes except on the coarsest level where an LU factorization is used.

In \cref{tab:signorini3d}, we test the effectiveness of the three-level OS-MG-PCG-S solver (block preconditioner \cref{fig:block} with the operator Schur complement approximant \cref{eq:So} and multigrid with PCG relaxation and vertex-star patch preconditioner) as $\lambda \to \infty$. We change from a Chebyshev relaxation to a PCG relaxation as (unreported) numerical experiments revealed that the Chebyshev relaxation was no longer effective at producing accurate multigrid cycles in this regime. We consider three discretizations for $(u,\psi)$: $(\mathrm{CG}_1)^3 \times \mathrm{CG}_1$ (denoted CG), $(\mathrm{CG}^A_3)^3 \times \mathrm{CG}_3$ (denoted CG-A), and $\mathrm{MTW}_1 \times \mathrm{DG}_0$ (denoted MTW). For the CG and MTW discretizations, we use the default prolongation operator, whereas a specialized prolongation is used for CG-A \cite[Sec.~3]{farrell2022}.

When the PG solve converged, it converged in six proximal steps or less. The average number of Newton iterations per proximal step was 2.60 or less with a CG and CG-A discretization, and 4.00 or less with an MTW discretization.

When $\lambda = 10, 10^2$ and $10^3$, the CG and CG-A discretizations produce similar average FGMRES iteration counts per Newton step, ranging from 5.15 to 8.15 with CG and 4.62 to 5.29 with CG-A. However, once $\lambda=10^4$ (Poisson's ratio of 0.49), the OS-MG-PCG-S solver failed to converge with the CG discretization. In contrast, we achieve up to $\lambda=10^5$ (Poisson's ratio of 0.499) with the CG-A discretization with an average of 8.38 FGMRES iterations per Newton step. The number of PCG iterations per $\hat{D}$-inverse was at most 9.28 across all choices for $\lambda$ with CG-A. The solver also reached $\lambda=10^5$ with the MTW discretization, although the iteration counts were higher ranging between 12.85 to 22.56 FGMRES iterations per Newton step. However, since a piecewise constant discretization is used for $\psi$, the $\hat{D}$-block is inverted with a single PCG iteration.

\subsection{Strain-constrained problem}
\label{sec:examples:strain}

In this section, we consider a realization of the strain-constrained problem \cref{eq:strain}.

\begin{figure}[h!]
\centering
\includegraphics[width =0.95\textwidth]{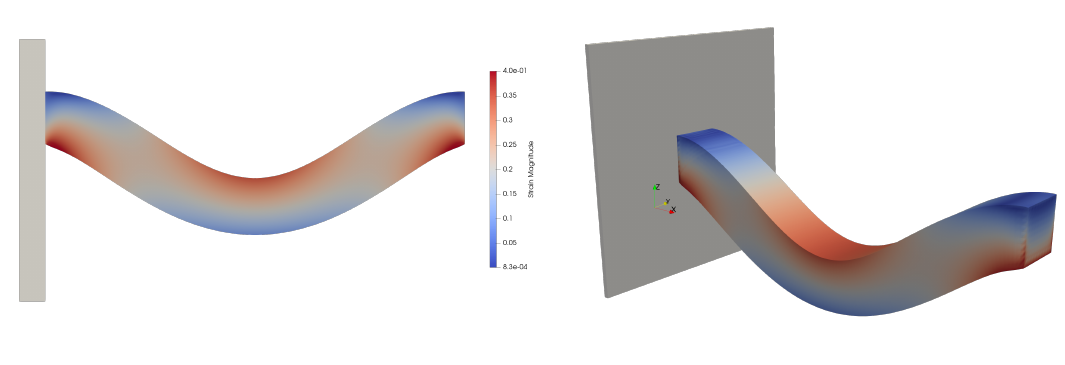}
\caption{The 2D (left) and 3D (right) deformed strain-constrained elastic beams described in \Cref{sec:examples:strain}. The discretization for $(u,\psi)$ is $(\mathrm{CG}_1)^d \times (\mathrm{DG}_0)^{d \times d}_{\mathrm{sym}}$.}
\label{fig:strain}
\end{figure}

In both 2D and 3D, the setup is
\begin{align}
\begin{gathered}
\Omega = (0,1) \times (0,0.1)^{d-1}, \quad f = -10\vectt{e}_d, \quad \mu = 70, \quad \lambda =100,\\
\phi = 0.4, \quad \Gamma_D = \{ x_1 = 0\} \cup \{x_1=1\}, \quad g = 
\begin{cases}
\vectt{0} & x_1 = 0,\\
-0.2 \vectt{e}_1 & x_1 = 1.
\end{cases}
\end{gathered}
\end{align}
Here $\vectt{e}_j \in \mathbb{R}^d$ denotes a unit vector with a 1 in the $j$-th entry and 0 elsewhere. This setup models a thin beam, with Poisson's ratio $\nu = 0.294$, fixed to a wall on the left-hand edge and face in two and three dimensions, respectively. The right-hand side edge (resp.~face) is compressed towards the wall by a value of 0.2 and the beam experiences a downwards body force which together cause the beam to deflect downwards. However, the magnitude of the strain cannot exceed the value of $\phi=0.4$ pointwise. We plot approximations of the solutions in \cref{fig:strain}.

For this example, we solely consider a continuous piecewise linear discretization for the vector-valued function $u$ and a piecewise constant discretization for the symmetric tensor-valued function $\psi$, i.e.~$(\mathrm{CG}_1)^d \times (\mathrm{DG}_0)^{d \times d}_{\mathrm{sym}}$. In (unreported) experiments, we found that higher-order discretizations were either unstable or converged to results that significantly violated the constraints. As the focus of this paper is to develop preconditioners for the PG linear systems, we leave the investigation of additional stable pairs for the strain-constrained problem for later work. 

All runs are initialized with $\alpha_0 = 10^{-2}$ in 2D and $\alpha_0=10^{-3}$ in 3D. We use the proximal update rule $\alpha_k = \min\{2 \alpha_{k-1}, 10\}$, and terminate once $\|u^k-u^{k-1}\|_{H^1(\Omega)} \leq 10^{-4}$. We pick the parameter $\gamma$ in \cref{eq:pg:strain} as $\gamma=100$ in 2D and $\gamma=5000$ in 3D. We pick $u^0 = \psi^0 = 0$  where the Dirichlet boundary conditions are pre-applied to $u^0$ before the first Newton solve. The initial iterate for Newton's method at the $k^\text{th}$-proximal step is $(u^{k-1}, \psi^{k-1})$.  The Newton solves terminate when the $\ell^2$-norm of the residual drops below $10^{-6}$ and we use an absolute and relative tolerance of $10^{-7}$ for the Krylov solvers. We apply five relaxation iterations on each level of the MG schemes except on the coarsest level where an LU factorization is used.  

In \cref{tab:strain}, on a base mesh and its uniform refinement, we report the results of an LU factorization in 2D as well as the block preconditioner \cref{fig:block} with either the algebraic \cref{eq:Sd} or operator \cref{eq:So} Schur complement approximant inverted by a multigrid cycle with GMRES relaxation and a block-Jacobi preconditioner for both 2D and 3D. Due to the severe nonlinearity of the problem and the fact that $\psi$ is now tensor-valued, we found that other relaxation iterates such as Chebyshev or PCG no longer lead to robust multigrid cycles. We also tested a monolithic Vanka-type strategy, but this also did not converge.

\subsubsection{2D Strain} All PG solves converged in 12 proximal steps. Moreover, the number of Newton steps per proximal step remained consistent across the three solvers, namely 2.08 Newton iterations per proximal step on the base mesh and 1.92 Newton iterations on the uniform refinement. The average number of FGMRES iterations per Newton step was reasonable, under seven for both block preconditioners. The operator Schur complement approximant produced slightly better iteration counts. Since the discretization for $\psi$ is piecewise constant, the solver only required one PCG iteration for each $\hat{D}$-block inverse.

\subsubsection{3D Strain} In 3D, the PG solves converged in 17 proximal steps with an average of at most 2.06 Newton iterations per proximal step on both meshes. The average number of FGMRES iterations per Newton step was 5.20 or less, with the operator Schur complement approximant resulting in slightly smaller iteration counts than the algebraic Schur complement approximant.

\FloatBarrier

\begin{table}[h!]
\renewcommand{\arraystretch}{1.2}
\centering
\scalebox{0.75}{
\begin{tabular}{|l|c|c|c|c|c|c|c|}
\hhline{|-|-|-|-|-|-|-|-|}
\rowcolor{lightgray!10}
\multicolumn{1}{c|}{LS Solver} & \multicolumn{1}{c|}{$N$} & \multicolumn{1}{c|}{$p$} & \multicolumn{1}{c|}{DoFs} & \multicolumn{1}{c|}{Levels} & \multicolumn{1}{c|}{PG} & \multicolumn{1}{c|}{Newton} & \multicolumn{1}{c|}{FGMRES} \\
\hline
\multirow{9}{*}{LU} & 64 & 1 & 31,158 & - & 9 & 1.78 & -\\
& 64 & 2 & 123,546 & - & 9 & 1.67 & -  \\
 & 64 & 3 & 277,166 & - & 9 & 1.67 & -  \\
 & 128 & 1 & 124,354 & -& 9 & 1.67 & -  \\
 & 128 & 2 & 495,250 & - & 9 & 1.67 & -  \\
 & 128 & 3 & 1,112,690 & - & 9 & 1.56 & - \\
 & 256 & 1 & 495,250 & - & 9 & 1.56 & - \\
 & 256 & 2 & 1,976,674 & - & 9 & 1.56 & - \\
 & 256 & 3 & 4,444,274 & - & 9 & 1.44 & -  \\
 \hline
\multirow{12}{*}{MG-Vanka}  & 64 & 1 & 31,158 & 2 & 9 & 1.78 & 2.94 (4)  \\
& 64 & 2 & 123,546 & 2 & 9 & 1.67 & 2.53 (4)  \\
& 64 & 3 & 277,166 & 2 & 9 & 1.67 & 2.33 (3)  \\
 & 128 & 1 & 123,546 & 3 & 9 & 1.67 & 2.87 (3) \\
 & 128 & 2 & 492,018 & 3 & 9 & 1.67 & 2.47 (4)  \\
 & 128 & 3 & 1,105,418 & 3 & 9 & 1.56 & 2.36 (3) \\
& 128 & 1 & 124,354 & 2 & 9 & 1.67 & 2.73 (3)  \\
 & 128 & 2 & 495,250 & 2 & 9 & 1.67 & 2.40 (3) \\
 & 128 & 3 & 1,112,690 & 2 & 9 & 1.56 & 2.36 (3)  \\
& 256 & 1 & 495,250 & 3 & 9 & 1.56 & 2.71 (3) \\
& 256 & 2 & 1,976,674 & 3 & 9 & 1.56 & 2.43 (3) \\
 & 256 & 3 & 4,444,274 & 3 & 9 & 1.44 & 2.31 (3)  \\
\hline
\end{tabular}}
\caption{(Obstacle Problem, 2D). A comparison of LU and the monolithic multigrid solver  (\cref{fig:vanka})  for the PG systems of the two-dimensional obstacle problem in \Cref{sec:examples:obstacle}. Here $\epsilon = 0$ and a $\mathrm{CG}_p \times \mathrm{CG}_p$ discretization is used for $(u,\psi)$ with $p \in \{1,2,3\}$.}
\label{tab:obstacle2d-1}
\end{table}

\begin{table}[!h]
\renewcommand{\arraystretch}{1.2}
\centering
\scalebox{0.78}{
\begin{tabular}{|l|c|c|c|c|c|c|c|c|}
\hhline{|-|-|-|-|-|-|-|-|-|}
\rowcolor{lightgray!10}
\multicolumn{1}{c|}{LS Solver} & \multicolumn{1}{c|}{$N$} & \multicolumn{1}{c|}{$p$} & \multicolumn{1}{c|}{DoFs} & \multicolumn{1}{c|}{Levels} & \multicolumn{1}{c|}{PG} & \multicolumn{1}{c|}{Newton} & \multicolumn{1}{c|}{FGMRES} & \multicolumn{1}{c|}{PCG} \\
\hline
\multirow{12}{*}{AS-MG-Ch-J}& 64 & 1 & 31,158 & 2 & 9 & 1.78 & 4.94 (9) & 3.97 (4) \\
& 64 & 2 & 123,546 & 2 & 9 & 1.67 & 4.13 (9) & 5.08 (6) \\
& 64 & 3 & 277,166 & 2 & 9 & 1.67 & 3.87 (9) & 5.34 (7) \\
& 128 & 1 & 123,546 & 3 & 9 & 1.67 & 3.93 (8) & 3.97 (4) \\
& 128 & 2 & 492,018 & 3 & 9 & 1.67 & 3.67 (8) & 4.95 (5) \\
& 128 & 3 & 1,105,418 & 3 & 9 & 1.56 & 3.50 (7) & 4.98 (6) \\
 & 128 & 1 & 124,354 & 2 & 9 & 1.67 & 3.93 (8) & 3.97 (4) \\
& 128 & 2 & 495,250 & 2 & 9 & 1.67 & 3.53 (8) & 4.95 (5) \\
& 128 & 3 & 1,112,690 & 2 & 9 & 1.56 & 3.36 (7) & 5.00 (6) \\
& 256 & 1 & 495,250 & 3 & 9 & 1.56 & 3.57 (6) & 3.96 (4) \\
& 256 & 2 & 1,976,674 & 3 & 9 & 1.56 & 3.50 (6) & 4.95 (5) \\
& 256 & 3 & 4,444,274 & 3 & 9 & 1.44 & 3.46 (6) & 4.94 (5) \\
\hline
\multirow{12}{*}{AS-MG-Ch-S} & 64 & 1 & 31,158 & 2 & 9 & 1.78 & 4.88 (9) & 3.97 (4) \\
 & 64 & 2 & 123,546 & 2 & 9 & 1.67 & 3.93 (9) & 5.08 (6) \\
 & 64 & 3 & 277,166 & 2 & 9 & 1.67 & 3.67 (8) & 5.32 (7) \\
& 128 & 1 & 123,546 & 3 & 9 & 1.67 & 3.93 (8) & 3.97 (4) \\
& 128 & 2 & 492,018 & 3 & 9 & 1.67 & 3.33 (7) & 4.95 (5) \\
& 128 & 3 & 1,105,418 & 3 & 9 & 1.56 & 3.14 (5) & 4.94 (6) \\
 & 128 & 1 & 124,354 & 2 & 9 & 1.67 & 3.73 (8) & 3.96 (4) \\
& 128 & 2 & 495,250 & 2 & 9 & 1.67 & 3.20 (6) & 4.95 (5) \\
& 128 & 3 & 1,112,690 & 2 & 9 & 1.56 & 3.14 (5) & 4.97 (6) \\
 & 256 & 1 & 495,250 & 3 & 9 & 1.56 & 3.29 (6) & 3.96 (4) \\
& 256 & 2 & 1,976,674 & 3 & 9 & 1.56 & 2.71 (4) & 4.93 (5) \\
& 256 & 3 & 4,444,274 & 3 & 9 & 1.44 & 2.69 (4) & 4.93 (5) \\
\hline
\multirow{12}{*}{OS-MG-Ch-J}& 64 & 1 & 31,158 & 2 & 9 & 1.78 & 4.62 (9) & 3.97 (4) \\
& 64 & 2 & 123,546 & 2 & 9 & 1.67 & 3.67 (10) & 5.11 (6) \\
& 64 & 3 & 277,166 & 2 & 9 & 1.67 & 3.33 (9) & 5.37 (7) \\
& 128 & 1 & 123,546 & 3 & 9 & 1.67 & 3.33 (8) & 3.96 (4) \\
& 128 & 2 & 492,018 & 3 & 9 & 1.67 & 3.00 (8) & 4.94 (5) \\
& 128 & 3 & 1,105,418 & 3 & 9 & 1.56 & 2.50 (6) & 4.96 (6) \\
& 128 & 1 & 124,354 & 2 & 9 & 1.67 & 3.40 (8) & 3.96 (4) \\
& 128 & 2 & 495,250 & 2 & 9 & 1.67 & 3.00 (8) & 4.94 (5) \\
 & 128 & 3 & 1,112,690 & 2 & 9 & 1.56 & 2.43 (6) & 5.00 (6) \\
 & 256 & 1 & 495,250 & 3 & 9 & 1.56 & 2.86 (6) & 3.95 (4) \\
& 256 & 2 & 1,976,674 & 3 & 9 & 1.56 & 2.21 (4) & 4.92 (5) \\
& 256 & 3 & 4,444,274 & 3 & 9 & 1.44 & 1.77 (3) & 4.89 (5) \\
\hline
\multirow{12}{*}{OS-MG-Ch-S}& 64 & 1 & 31,158 & 2 & 9 & 1.78 & 4.62 (9) & 3.97 (4) \\
 & 64 & 2 & 123,546 & 2 & 9 & 1.67 & 3.67 (10) & 5.11 (6) \\
 & 64 & 3 & 277,166 & 2 & 9 & 1.67 & 3.33 (9) & 5.38 (7) \\
 & 128 & 1 & 123,546 & 3 & 9 & 1.67 & 3.33 (8) & 3.96 (4) \\
 & 128 & 2 & 492,018 & 3 & 9 & 1.67 & 3.00 (8) & 4.94 (5) \\
 & 128 & 3 & 1,105,418 & 3 & 9 & 1.56 & 2.57 (6) & 4.96 (6) \\
& 128 & 1 & 124,354 & 2 & 9 & 1.67 & 3.40 (8) & 3.96 (4) \\
 & 128 & 2 & 495,250 & 2 & 9 & 1.67 & 3.00 (8) & 4.94 (5) \\
& 128 & 3 & 1,112,690 & 2 & 9 & 1.56 & 2.57 (6) & 4.97 (6) \\
& 256 & 1 & 495,250 & 3 & 9 & 1.56 & 2.86 (6) & 3.95 (4) \\
 & 256 & 2 & 1,976,674 & 3 & 9 & 1.56 & 2.07 (4) & 4.91 (5) \\
 & 256 & 3 & 4,444,274 & 3 & 9 & 1.44 & 1.85 (3) & 4.90 (5) \\
\hline
\end{tabular}}
\caption{(Obstacle Problem, 2D). A comparison of block preconditioners from \cref{fig:block} (see \cref{tab:solvers}) for the two-dimensional obstacle problem in \Cref{sec:examples:obstacle} with an equal-order $\mathrm{CG}_p \times \mathrm{CG}_p$ discretization with $p \in \{1,2,3\}$. Here $\epsilon = 10^{-5}$.}
\label{tab:obstacle2d-2}
\end{table}

\begin{table}[!h]
\renewcommand{\arraystretch}{1.2}
\centering
\scalebox{0.84}{
\begin{tabular}{|l|c|c|c|c|c|c|c|c|c|}
\hhline{|-|-|-|-|-|-|-|-|-|-|}
\rowcolor{lightgray!10}
\multicolumn{1}{c|}{LS Solver} & \multicolumn{1}{c|}{$N$} & \multicolumn{1}{c|}{$p$} & \multicolumn{1}{c|}{DoFs} & \multicolumn{1}{c|}{$\epsilon$} & \multicolumn{1}{c|}{Levels} & \multicolumn{1}{c|}{PG} & \multicolumn{1}{c|}{Newton} & \multicolumn{1}{c|}{FGMRES} & \multicolumn{1}{c|}{PCG} \\
\hline
\multirow{4}{*}{MG-Vanka}
& 40 & 1 & 702,746 & 0 & 2 & 25 & 1.60 & 2.52 (3) & - \\
& 40 & 2 & 5,441,530 & 0 & 2 & 25 & 1.52 & 2.18 (3) & - \\
& 80 & 1 & 5,441,530 & 0 & 3 & 25 & 1.48 & 2.62 (4) & - \\
& 80 & 2 & 42,813,778 & 0 & 3 & 25 & 1.36 & 2.26 (3) & - \\
\hline
\multirow{4}{*}{OS-MG-Ch-J}
& 40 & 1 & 702,746 & 1e-4 & 2 & 25 & 1.60 & 3.75 (9) & 7.70 (8) \\
& 40 & 2 & 5,441,530 & 1e-4 & 2 & 25 & 1.52 & 3.00 (7) & 14.25 (19) \\
& 80 & 1 & 5,441,530 & 1e-4 & 3 & 25 & 1.48 & 3.03 (6) & 7.46 (8) \\
& 80 & 2 & 42,813,778 & 1e-4 & 3 & 25 & 1.36 & 2.53 (5) & 12.47 (14) \\
\hline
\multirow{3}{*}{OS-MG-Ch-S}
& 40 & 1 & 702,746 & 1e-4 & 2 & 25 & 1.60 & 3.75 (9) & 7.72 (8) \\
& 40 & 2 & 5,441,530 & 1e-4 & 2 & 25 & 1.52 & 2.87 (7) & 14.62 (20) \\
& 80 & 1 & 5,441,530 & 1e-4 & 3 & 25 & 1.48 & 3.03 (6) & 7.15 (8) \\
\hline
\multirow{2}{*}{AS-MG-Ch-J}
& 40 & 1 & 702,746 & 1e-4 & 2 & 25 & 1.60 & 4.05 (8) & 7.76 (8) \\
& 80 & 1 & 5,441,530 & 1e-4 & 3 & 25 & 1.48 & 3.59 (6) & 7.12 (8) \\
\hline
\multirow{2}{*}{AS-MG-Ch-S}
& 40 & 1 & 702,746 & 1e-4 & 2& 25 & 1.60 & 3.95 (8) & 7.78 (8) \\
& 80 & 1 & 5,441,530 & 1e-4 & 3 & 25 & 1.48 & 3.24 (6) & 7.07 (8) \\
\hline
\end{tabular}}
\caption{(Obstacle Problem, 3D). A comparison of iterative solvers (\cref{fig:vanka,fig:block})  for the three-dimensional obstacle problem in \Cref{sec:examples:obstacle} with an equal-order $\mathrm{CG}_p \times \mathrm{CG}_p$ discretization with $p \in \{1,2\}$.}
\label{tab:obstacle3d}
\end{table}

\begin{table}[!h]
\renewcommand{\arraystretch}{1.2}
\centering
\scalebox{0.78}{
\begin{tabular}{|l|c|c|c|c|c|c|c|c|c|c|}
\hhline{|-|-|-|-|-|-|-|-|-|-|-|}
\rowcolor{lightgray!10}
\multicolumn{1}{c|}{FEM} & \multicolumn{1}{c|}{LS Solver} & \multicolumn{1}{c|}{Levels} & \multicolumn{1}{c|}{$N$} & \multicolumn{1}{c|}{$p$} & \multicolumn{1}{c|}{$u$ DoFs} & \multicolumn{1}{c|}{$\psi$ DoFs} & \multicolumn{1}{c|}{PG} & \multicolumn{1}{c|}{Newton} & \multicolumn{1}{c|}{FGMRES} & \multicolumn{1}{c|}{PCG} \\
\hline
\multirow{4}{*}{CG} & \multirow{4}{*}{LU} & - & 60 & 1 & 34,438 & 189 & 6 & 2.67 & -& - \\
&  & -& 60 & 2 & 136,730 & 377 & 6 & 2.67 & - & - \\
 & & -& 120 & 1 & 136,730 & 377 & 6 & 2.67 & -& - \\
& & - & 120 & 2 & 544,882 & 753 & 6 & 2.50 & -& - \\
\hline
\multirow{4}{*}{CG} & \multirow{4}{*}{OS-MG-Ch-J} & 2 & 60 & 1 & 34,438 & 189 & 6 & 2.67 & 4.19 (9) & 1.70 (3) \\
&  & 2 & 60 & 2 & 136,730 & 377 & 6 & 2.67 & 4.62 (9) & 1.55 (3) \\
 & & 3 & 120 & 1 & 136,730 & 377 & 6 & 2.67 & 4.31 (9) & 1.69 (3) \\
& & 3 & 120 & 2 & 544,882 & 753 & 6 & 2.50 & 4.33 (9) & 1.55 (3) \\
\hline
CG &\multirow{8}{*}{OS-MG-Ch-S} & 2 & 60 & 1 & 34,438 & 189 & 6 & 2.67 & 4.06 (8) & 1.71 (3) \\
MTW & & 2 & 60 & 1 & 153,438 & 188 & 6 & 2.67 & 9.69 (14) & 1.00 (1) \\
CG & & 2 & 60 & 2 & 136,730 & 377 & 6 & 2.67 & 3.62 (8) & 1.41 (3) \\
CG-A & & 2 & 60 & 2 & 408,154 & 377 & 6 & 2.67 & 3.31 (7) & 1.42 (3) \\
CG & & 3 & 120 & 1 & 136,730 & 377 & 6 & 2.67 & 4.06 (9) & 1.65 (3) \\
MTW &  & 3 & 120 & 1 & 612,228 & 376 & 6 & 2.67 & 11.75 (22) & 1.00 (1) \\
CG & & 3 & 120 & 2 & 544,882 & 753 & 6 & 2.50 & 3.00 (7) & 1.46 (3) \\
CG-A &  & 3 & 120 & 2 & 1,630,578 & 753 & 6 & 2.50 & 2.73 (6) & 1.46 (3) \\
\hline
\end{tabular}}
\caption{(Signorini Problem, 2D). A comparison of LU and block preconditioners from \cref{fig:block}  for the two-dimensional Signorini problem in \Cref{sec:examples:signorini} where $\lambda=10$. CG, CG-A and MTW denote an equal-order $(\mathrm{CG}_p)^2 \times \mathrm{CG}_p$, Alfeld-split $(\mathrm{CG}^A_2)^2 \times \mathrm{CG}_2$, and $\mathrm{MTW}_1 \times \mathrm{DG}_0$ discretization for $(u,\psi)$, respectively. For the LU solver we pick $\epsilon=0$ and in all the other linear solvers  $\epsilon = 10^{-4}$.}
\label{tab:signorini2d-1}
\end{table}

\begin{table}[!h]
\renewcommand{\arraystretch}{1.2}
\centering
\scalebox{0.88}{
\begin{tabular}{|l|c|c|c|c|c|c|c|}
\hhline{|-|-|-|-|-|-|-|-|}
\rowcolor{lightgray!10}
\multicolumn{1}{c|}{FEM} & \multicolumn{1}{c|}{$\lambda$} & \multicolumn{1}{c|}{$u$ DoFs} & \multicolumn{1}{c|}{$\psi$ DoFs} & \multicolumn{1}{c|}{PG} & \multicolumn{1}{c|}{Newton} & \multicolumn{1}{c|}{FGMRES} & \multicolumn{1}{c|}{PCG} \\
\hline
\multirow{4}{*}{CG}& 10 & 138,406 & 377 & 6 & 2.67 & 3.88 (8) & 1.66 (3) \\
& 100 & 138,406 & 377 & 6 & 2.50 & 8.53 (19) & 1.77 (3) \\
& 1000 & 138,406 & 377 & 6 & 2.50 & 26.33 (53) & 1.82 (3) \\
& 10,000 & 138,406 & 377 & $\times$ &  $\times$ &  $\times$&  $\times$ \\
\hline
\multirow{5}{*}{CG-A}& 10 & 408,154 & 377 & 6 & 2.67 & 3.31 (7) & 1.42 (3) \\
& 100 & 408,154 & 377 & 6 & 2.50 & 3.87 (7) & 1.54 (3) \\
& 1000 & 408,154 & 377 & 6 & 2.50 & 4.47 (9) & 1.52 (3) \\
& 10,000 & 408,154 & 377 & 6 & 2.33 & 4.43 (8) & 1.52 (3) \\
& 100,000 & 408,154 & 377 & 6 & 2.17 & 4.92 (9) & 1.45 (3) \\
\hline
\multirow{5}{*}{MTW} & 10 & 153,438 & 188 & 6 & 2.67 & 9.69 (14) & 1.00 (1) \\
& 100 & 153,438 & 188 & 6 & 2.83 & 9.53 (14) & 1.00 (1) \\
& 1000 & 153,438 & 188 & 6 & 2.67 & 10.19 (14) & 1.00 (1) \\
& 10,000 & 153,438 & 188 & 6 & 2.50 & 10.13 (14) & 1.00 (1) \\
& 100,000 & 153,438 & 188 & 6 & 2.50 & 10.00 (14) & 1.00 (1) \\
\hline
\end{tabular}}
\caption{(Signorini Problem, 2D). An investigation into $\lambda$-robustness of the two-level OS-MG-Ch-S preconditioner from \cref{fig:block}  for the two-dimensional Signorini problem in \Cref{sec:examples:signorini}. CG, CG-A and MTW denote an equal-order $(\mathrm{CG}_1)^2 \times \mathrm{CG}_1$, Alfeld-split $(\mathrm{CG}^A_2)^2 \times \mathrm{CG}_2$, and $\mathrm{MTW}_1 \times \mathrm{DG}_0$ discretization for $(u,\psi)$, respectively. Here $\epsilon = 10^{-4}$, $N=60$ for CG-A and MTW and $N=120$ for CG. The iteration counts remain robust in the CG-A and MTW discretizations whereas they grow when using CG. The $\times$ denotes a failure to converge.}
\label{tab:signorini2d-2}
\end{table}

\begin{table}[!h]
\renewcommand{\arraystretch}{1.2}
\centering
\scalebox{0.88}{
\begin{tabular}{|l|c|c|c|c|c|c|c|}
\hhline{|-|-|-|-|-|-|-|-|}
\rowcolor{lightgray!10}
\multicolumn{1}{c|}{FEM} & \multicolumn{1}{c|}{$\lambda$} & \multicolumn{1}{c|}{$u$ DoFs} & \multicolumn{1}{c|}{$\psi$ DoFs} & \multicolumn{1}{c|}{PG} & \multicolumn{1}{c|}{Newton} & \multicolumn{1}{c|}{FGMRES} & \multicolumn{1}{c|}{PCG} \\
\hline
\multirow{4}{*}{CG}& 10 & 6,362,847 & 45,965 & 6 & 2.33 & 6.21 (23) & 6.28 (7) \\
& 100 & 6,362,847 & 45,965 & 5 & 2.60 & 5.15 (10) & 6.49 (7) \\
 & 1000 & 6,362,847 & 45,965 & 5 & 2.60 & 8.15 (18) & 6.38 (7) \\ 
& 10,000 &  6,362,847& 45,965  & $\times$ & $\times$ & $\times$ & $\times$ \\ \hline
\multirow{5}{*}{CG-A} & 10 & 10,170,453 & 25,027 & 6 & 2.33 & 4.79 (10) & 8.54 (11) \\
& 100 & 10,170,453 & 25,027 & 5 & 2.60 & 4.62 (9) & 9.28 (12) \\
& 1000 & 10,170,453 & 25,027 & 5 & 2.80 & 5.29 (10) & 8.20 (10) \\
 & 10,000 & 10,170,453 & 25,027 & 5 & 2.60 & 6.31 (11) & 8.73 (11) \\
& 100,000 & 10,170,453 & 25,027 & 5 & 2.60 & 8.38 (15) & 8.34 (11) \\ \hline
\multirow{5}{*}{MTW}& 10 & 2,108,544 & 5520 & 5 & 2.60 & 12.85 (18) & 1.00 (1) \\
& 100 & 2,108,544 & 5520 & 5 & 2.80 & 13.00 (21) & 1.00 (1) \\
& 1000 & 2,108,544 & 5520 & 4 & 3.75 & 13.87 (21) & 1.00 (1) \\
 & 10,000 & 2,108,544 & 5520 & 4 & 4.00 & 18.38 (47) & 1.00 (1) \\
& 100,000 & 2,108,544 & 5520 & 4 & 4.00 & 22.56 (53) & 1.00 (1) \\
\hline
\end{tabular}}
\caption{(Signorini Problem, 3D). An investigation into $\lambda$-robustness of the three-level OS-MG-PCG-S preconditioner (see \cref{tab:solvers}) for the three-dimensional Signorini problem in \Cref{sec:examples:signorini}. CG, CG-A and MTW denote an equal-order $(\mathrm{CG}_1)^3 \times \mathrm{CG}_1$, Alfeld-split $(\mathrm{CG}^A_3)^3 \times \mathrm{CG}_3$, and $\mathrm{MTW}_1 \times \mathrm{DG}_0$ discretization for $(u,\psi)$, respectively. Here $\epsilon = 10^{-4}$, $N=5$ for CG-A and MTW and $N=20$ for CG. The $\times$ denotes a failure to converge.}
\label{tab:signorini3d}
\end{table}

\clearpage
\begin{table}[!h]
\renewcommand{\arraystretch}{1.2}
\centering
\scalebox{0.8}{
\begin{tabular}{|l|c|c|c|c|c|c|c|c|}
\hhline{|-|-|-|-|-|-|-|-|-|}
\rowcolor{lightgray!10}
\multicolumn{1}{c|}{LS Solver} & \multicolumn{1}{c|}{$N$} & \multicolumn{1}{c|}{$u$ DoFs} & \multicolumn{1}{c|}{$\psi$ DoFs} & \multicolumn{1}{c|}{Levels} & \multicolumn{1}{c|}{PG} & \multicolumn{1}{c|}{Newton} & \multicolumn{1}{c|}{FGMRES} & \multicolumn{1}{c|}{PCG} \\
\hline
\multicolumn{9}{|c|}{2D}\\
\hline
\multirow{2}{*}{LU} & 80 & 129,762 & 384,000 & -  & 12 & 2.08 & - & - \\
 & 160 & 515,522 & 1,536,000 & - & 12 & 1.92 & -& - \\
\hline
\multirow{2}{*}{OS-MG-GMRES-bJ} & 80 & 129,762 & 384,000 & 2  & 12 & 2.08 & 4.32 (12) & 1 (1) \\
& 160 & 515,522 & 1,536,000 & 3 & 12 & 1.92 & 5.78 (19) & 1 (1) \\
\hline
\multirow{2}{*}{AS-MG-GMRES-bJ}  & 80 & 129,762 & 384,000 & 2  & 12 & 2.08 & 5.00 (14) & 1 (1) \\
& 160 & 515,522 & 1,536,000 & 3  & 12 & 1.92 & 6.57 (22) & 1 (1) \\
\hline
\multicolumn{9}{|c|}{3D}\\
\hline
\multirow{2}{*}{OS-MG-GMRES-bJ} & 20 & 265,923 & 2,880,000 & 2 & 17 & 2.06 & 5.00 (24) & 1 (1) \\
& 40 & 2,022,243 & 23,040,000 & 3 & 17 & 2.06 & 4.43 (21) & 1 (1) \\
\hline
\multirow{2}{*}{AS-MG-GMRES-bJ} & 20 & 265,923 & 2,880,000 & 2 & 17 & 2.06 & 5.20 (25) & 1 (1) \\
& 40 & 2,022,243 & 23,040,000 & 3 & 17 & 2.00 & 4.74 (22) & 1 (1) \\
\hline
\end{tabular}}
\caption{(Strain Problem, 2D \& 3D). A comparison of LU and block preconditioners from \cref{fig:block}  for the two- and three-dimensional strain-constrained problem in \Cref{sec:examples:strain}. We solely consider a $(\mathrm{CG}_1)^d \times (\mathrm{DG}_0)^{d \times d}_{\text{sym}}$, $d \in \{2,3\}$, discretization for $(u,\psi)$. For $\gamma$ in \cref{eq:pg:strain}, we choose $\gamma=100$ in 2D and $\gamma=5000$ in 3D. For the LU factorization we pick $\epsilon=0$ and in the other linear solvers  $\epsilon = 10^{-3}$.}
\label{tab:strain}
\end{table}

\FloatBarrier
\section{Conclusions}
\label{sec:conclusions}

In this paper, we investigated iterative Krylov methods for the Newton linear systems that arise when the proximal Galerkin solver is applied to the obstacle problem, Signorini problem, and strain-constrained elasticity problem. 

The first step is to utilize an FGMRES Krylov method to invert the Jacobian matrix in the Newton linear system. For FGMRES to converge within a reasonable number of iterations, we then construct a block preconditioner. Rather unusually, we perform a Schur complement factorization with respect to the degenerate latent block $\hat{D}$ (analogous to the pressure block in Stokes). This requires two inverses of the degenerate latent block and one of the primal Schur complement per FGMRES iteration. To stabilize the inverse of the degenerate latent block, we augment the block with a small mass-matrix perturbation. Note this perturbation \emph{is only included in the preconditioner}. The action of the original Newton matrix in the FGMRES iteration is left unperturbed and hence \emph{a true Newton step is calculated}. The stabilized degenerate block is analogous to a weighted mass matrix and we found that it could be inverted efficiently with a conjugate gradient method preconditioned by a block-Jacobi iteration.

We then proposed two sparse approximations of the  Schur complement.  One approximant breaks the continuity in the finite element basis of the latent variable and directly assembles the corresponding \emph{algebraic} Schur complement approximant. Up to the latent block perturbation, this method is exact when the discretization for the latent variable is discontinuous. The other approximant examines the individual components of the Schur complement and designs a continuous PDE approximant which is then discretized. We call this the \emph{operator} Schur complement approximant. 

Both these Schur complement approximants are inverted via geometric multigrid. Depending on the problems, we consider Chebyshev, conjugate gradient, and GMRES relaxation schemes preconditioned by Jacobi, block-Jacobi and vertex-star patch iterations. In practice, the operator Schur complement approximant was both cheaper and led to smaller iteration counts. Similarly, the cost of Chebyshev relaxation with a Jacobi preconditioner was the cheapest. However, as demonstrated in \Cref{sec:examples} it does not always lead to  parameter-robust multigrid cycles. We also introduce a monolithic multigrid approach, with Vanka-type relaxation, of the full Jacobian.  This approach proved very effective for the obstacle problem but was less effective for the strain-constrained elasticity problem and thus we conclude that the block preconditioner approach is robust for a broader class of problems.

Experimentally, we observed that the block preconditioning strategy results in modest and robust iteration counts for different mesh sizes, discretization degrees and finite element families. Moreover, we saw that the block preconditioner strategy with vertex-star patch relaxation, coupled with specific discretizations, leads to moderate iteration counts even for nearly-incompressible elastic materials. This was exemplified when we achieved convergence for a three-dimensional Signorini problem, with Poisson's ratio of 0.499, discretized with over ten million DoFs for the displacement variable. The solver only required an average of 8.38 FGMRES iterations (and thus 8.38 multigrid cycles) per Newton step.

We remark that because this work introduces several linear solvers, a convergence analysis is not conducted. However, we hope that the promising experiments will encourage an analysis for the individual solvers in future work.

\bibliographystyle{siamplain}
\bibliography{references}

\end{document}